# A POLYMER IN A MULTI-INTERFACE MEDIUM

By Francesco Caravenna and Nicolas Pétrélis

*University of Padova and EURANDOM*

We consider a model for a polymer chain interacting with a sequence of equispaced flat interfaces through a pinning potential. The intensity $\delta \in \mathbb{R}$ of the pinning interaction is constant, while the interface spacing $T = T_N$ is allowed to vary with the size $N$ of the polymer. Our main result is the explicit determination of the scaling behavior of the model in the large $N$ limit, as a function of $(T_N)_N$ and for fixed $\delta > 0$. In particular, we show that a transition occurs at $T_N = O(\log N)$. Our approach is based on renewal theory.

## 1. Introduction and main results.

1.1. *The model.* In this paper, we study a $(1+1)$-dimensional model for a polymer chain dipped in a medium constituted by infinitely many horizontal interfaces. The possible configurations of the polymer are modeled by the trajectories $\{(i, S_i)\}_{i \geq 0}$ of the simple symmetric random walk on $\mathbb{Z}$, with law denoted by **P**, that is, $S_0 = 0$ and $(S_i - S_{i-1})_{i \geq 1}$ is an i.i.d. sequence of Bernoulli trials satisfying $\mathbf{P}(S_1 = \pm 1) = 1/2$. We assume that the interfaces are equispaced, that is, at the same distance $T \in 2\mathbb{N}$ from each other (note that $T$ is assumed to be even for notational convenience, due to the periodicity of the simple random walk).

The interaction between the polymer and the medium is described by the following Hamiltonian:

$$(1.1) \qquad H_{N,\delta}^T(S) := \delta \sum_{i=1}^N \mathbf{1}_{\{S_i \in T\mathbb{Z}\}} = \delta \sum_{k \in \mathbb{Z}} \sum_{i=1}^N \mathbf{1}_{\{S_i = kT\}},$$

where $N \in \mathbb{N}$ is the size of the polymer and $\delta \in \mathbb{R}$ is the intensity of the energetic reward (if $\delta > 0$) or penalty (if $\delta < 0$) that the polymer receives when touching the interfaces. More precisely, the model is defined by the following probability law $\mathbf{P}_{N,\delta}^T$ on $\mathbb{R}^{\mathbb{N}\cup\{0\}}$:

$$(1.2) \qquad \frac{d\mathbf{P}_{N,\delta}^T}{d\mathbf{P}}(S) := \frac{\exp(H_{N,\delta}^T(S))}{Z_{N,\delta}^T},$$

where $Z_{N,\delta}^T = \mathbf{E}(\exp(H_{N,\delta}^T(S)))$ is the normalizing constant, called the *partition function*.

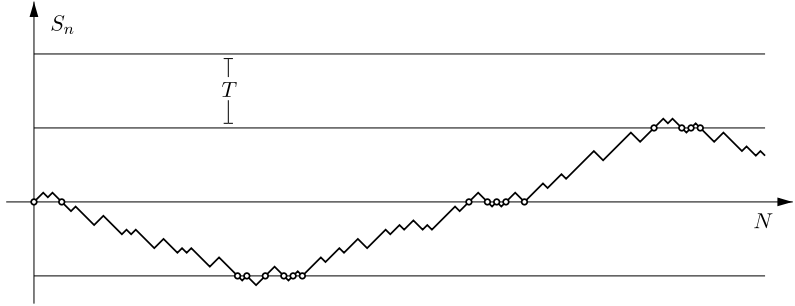

Fig. 1. *A typical path of the polymer measure $\mathbf{P}_{N,\delta}^T$ with $N = 158$ and $T = 16$. The circles represent the points where the polymer touches the interfaces, which are favored (resp., disfavored) when $\delta > 0$ (resp., $\delta < 0$).*

It should be clear that the effect of the Hamiltonian $H_{N,\delta}^T$ is to favor or penalize, according to the sign of $\delta$, the trajectories $\{(n, S_n)\}_n$ that have large numbers of intersections with the interfaces, located at heights $T\mathbb{Z}$ (we refer to Figure 1 for a graphical description). Although we present a number of results in this work that do not depend on the sign of $\delta$, we stress that, hereafter, our main concern is with the case $\delta > 0$.

If we let $T \to \infty$ in (1.2) *for fixed $N$* (in fact, it suffices to take $T > N$), we obtain a well-defined limiting model $\mathbf{P}_{N,\delta}^\infty$:

$$(1.3) \quad \frac{d\mathbf{P}_{N,\delta}^\infty}{d\mathbf{P}}(S) := \frac{\exp(H_{N,\delta}^\infty(S))}{Z_{N,\delta}^\infty}, \qquad \text{where } H_{N,\delta}^\infty(S) := \delta \sum_{i=1}^N \mathbf{1}_{\{S_i=0\}}.$$

$\mathbf{P}_{N,\delta}^\infty$ is known in the literature as a *homogeneous pinning model* and it describes a polymer chain interacting with a single flat interface, namely the $x$-axis. This model, together with several variants (such as the *wetting model*, where $\{S_n\}_n$ is also constrained to stay nonnegative), has been studied in depth, first in the physical literature (see [11] and references therein) and, more recently, in the mathematical literature [4, 9, 12, 13]. In particular, it is well known that a *phase transition* between a delocalized regime and a localized one occurs as $\delta$ varies and this transition can be characterized in terms of the path properties of $\mathbf{P}_{N,\delta}^\infty$.

The aim of this paper is to answer the same kinds of questions for the model $\mathbf{P}_{N,\delta}^{T_N}$, as a function of $\delta$ and of the interface spacing $T = T_N$, which is allowed to vary with $N$. We denote the full sequence by $\mathbf{T} := (T_N)_{N \in \mathbb{N}}$ (taking values in $2\mathbb{N}$) and, with no essential loss of generality (one could focus on subsequences), we assume that $\mathbf{T}$ has a limit as $N \to \infty$:

(1.4) $$\exists \lim_{N \to \infty} T_N =: T_\infty \in 2\mathbb{N} \cup \{+\infty\}.$$

[Of course, if $T_\infty < \infty$, then the sequence $(T_N)_N$ must eventually take the constant value $T_\infty$.] For notational convenience, we also assume that $T_N \leq N$: again, there is no real loss of generality since, for $T_N > N$, the law $\mathbf{P}_{N,\delta}^{T_N}$ reduces to the just-mentioned $\mathbf{P}_{N,\delta}^\infty$.

Before precisely stating the results we obtain in this paper, let us briefly describe the motivations behind our model and its context. Several models for a polymer interacting with a single linear interface have been investigated in the past twenty years, both in the physical and in the mathematical literature (see [11] and [12] for two excellent surveys). Probably the two most popular classes among them are the so-called *copolymer at a selective interface* separating two selective solvents and the *pinning of a polymer at an interface*, of which the homogeneous pinning model $\mathbf{P}_{N,\delta}^\infty$ is the simplest and most basic example. Although some questions still remain open, notably when *disorder* is present, important progress has been made and there is now a fairly good understanding of the mechanism leading to phase transitions for these models.

More recently, some generalizations have been introduced to account for interactions taking place on structures more general than a single linear interface. In the *copolymer class*, we mention [6, 7] and [17], where the medium is constituted by an emulsion, and, especially, [8], where the single linear interface is replaced by infinitely many equispaced flat interfaces, separating alternate layers of each selective solvent. Our model $\mathbf{P}_{N,\delta}^T$ provides a closely analogous generalization in the *pinning class*, with the important difference that the model considered in [8] is disordered. In a sense, what we consider is the simplest case of a pinning model interacting with infinitely many interfaces. In analogy with the single interface case [12], we believe that understanding this basic example in detail is the first step toward a comprehension of the more sophisticated disordered case.

Let us briefly describe the results obtained in [8]. The authors focus on the case where the interface spacing $T_N$ diverges as $N \to \infty$ and they show that the free energy of the model is the same as in the case of one single linear interface. Then, under stronger assumption on $(T_N)_N$, namely $T_N/\log\log N \to \infty$ and $T_N/\log N \to 0$, they show that the polymer visits infinitely many different interfaces and the asymptotic behavior of the time needed to hop from an interface to a neighboring behaves like $e^{cT_N}$.

In this paper, we consider analogous questions for our model $\mathbf{P}_{N,\delta}^T$. In our nondisordered setting, we obtain stronger results: in particular, we are able to precisely describe the path behavior of the polymer in the large $N$ limit for an arbitrary sequence $\mathbf{T} = (T_N)_N$ and for $\delta > 0$ (i.e., we focus on the case of attractive interfaces). In fact, there is a subtle interplay between the pinning reward $\delta$ and the speed $T_N$ at which the interfaces depart, which is responsible for the scaling behavior of the polymer. It turns out that there are three different regimes, determined by comparing $T_N$ with $\frac{\log N}{c_\delta}$, where $c_\delta > 0$ is computed explicitly. We refer to Theorem 2 and to the following discussion for a detailed explanation of our results. Let us just mention that, as $T_N$ increases from $O(1)$ to the critical speed $\frac{\log N}{c_\delta}$, the scaling constants of $S_N$ decrease smoothly from the diffusive behavior $\sqrt{N}$ to $\log N$, while if $T_N \gg \frac{\log N}{c_\delta}$, then $S_N = O(1)$. This means, on the one hand, that by accelerating the growth of the interface spacing, the scaling of $S_N$ decreases and, on the other hand, that scaling behaviors for $S_N$ intermediate between $O(1)$ and $\log N$ (such as, e.g., $\log \log N$) are not possible in our model. We also stress that our model is *subdiffusive* as soon as $T_N \to \infty$. Subdiffusive behaviors appear in a variety of models dealing with random walks subject to some form of penalization: from the (very rich) literature, we mention, for instance, [1] and [15] on the mathematical side and [16] on the physical side.

Our approach is mainly based on renewal theory. The use of these kinds of techniques in the field of polymer models has proven to be extremely successful, starting from [4] and [9], and, more recently, has been generalized to cover Markovian settings; see [2] and [5]. The key point is to get sharp estimates on suitable renewal functions.

The same approach can be applied to deal with the depinning case $\delta < 0$, that is, when touching an interface entails a penalty rather than a reward. However, in this case, the limiting model $\mathbf{P}_{N,\delta}^\infty$ is *delocalized* and this fact generates additional nontrivial difficulties. For this reason, the analysis of the $\delta < 0$ case is given in a separate paper [3], where we show that there are remarkable differences with respect to the $\delta > 0$ case that we consider here. In particular, the critical speed of $T_N$ above which the polymer gives up visiting infinitely many different interfaces is no longer of order $\log N$, but rather of order $N^{1/3}$ (see Theorem 1 in [3] for a precise statement).

1.2. *The free energy.* The standard way of studying the effect of the interaction (1.1) for large $N$ is to look at the *free energy* of the model, defined as the limit

$$(1.5) \quad \phi(\delta, \mathbf{T}) := \lim_{N \to \infty} \phi_N(\delta, \mathbf{T}), \qquad \text{where } \phi_N(\delta, \mathbf{T}) := \frac{1}{N} \log Z_{N,\delta}^{T_N}.$$

The existence of such a limit, for any choice of $\delta \in \mathbb{R}$ and $\mathbf{T}$ satisfying (1.4), is proved in Section 2. To understand why one should look at $\phi$, we introduce the random variable

$$(1.6) \qquad L_{N,T} := \sum_{i=1}^{N} \mathbf{1}_{\{S_i \in T\mathbb{Z}\}} = \sum_{k \in \mathbb{Z}} \sum_{i=1}^{N} \mathbf{1}_{\{S_i = kT\}}$$

and observe that an easy computation yields

$$\frac{\partial}{\partial \delta} \phi_N(\delta, \mathbf{T}) = \mathbf{E}_{N,\delta}^{T_N}\left(\frac{L_{N,T_N}}{N}\right),$$

$$\frac{\partial^2}{\partial \delta^2} \phi_N(\delta, \mathbf{T}) = N \operatorname{var}_{\mathbf{P}_{N,\delta}^{T_N}}\left(\frac{L_{N,T_N}}{N}\right) \geq 0.$$

In particular, $\phi_N(\delta, \mathbf{T})$ is a convex function of $\delta$ for every $N \in \mathbb{N}$. Hence, $\phi(\delta, \mathbf{T})$ is also convex and, by elementary convex analysis, it follows that provide $\phi(\delta, \mathbf{T})$ is differentiable,

$$(1.7) \qquad \frac{\partial}{\partial \delta} \phi(\delta, \mathbf{T}) = \lim_{N \to \infty} \mathbf{E}_{N,\delta}^{T_N}\left(\frac{L_{N,T_N}}{N}\right).$$

Thus, the first derivative of $\phi(\delta, \mathbf{T})$ gives the asymptotic proportion of time spent by the polymer on the interfaces, which explains the reason for looking at $\phi(\delta, \mathbf{T})$. In fact, a basic problem is the determination of the set of values of $\delta$ (if any) where $\phi(\delta, \mathbf{T})$ is not analytic, which correspond physically to the occurrence of a *phase transition* in the system.

This issue is addressed by our first result, which provides an explicit formula for $\phi(\delta, \mathbf{T})$. Let us introduce, for $T \in 2\mathbb{N} \cup \{+\infty\}$, the random variable $\tau_1^T$ defined by

$$(1.8) \qquad \tau_1^T := \inf\{n > 0 : S_n \in \{-T, 0, +T\}\}$$

and denote by $Q_T(\lambda)$ its Laplace transform under the simple random walk law $\mathbf{P}$:

$$(1.9) \qquad Q_T(\lambda) := \mathbf{E}(e^{-\lambda \tau_1^T}) = \sum_{n=1}^{\infty} e^{-\lambda n} \mathbf{P}(\tau_1^T = n).$$

When $T = +\infty$, the variable $\tau_1^\infty$ is nothing but the first return time of the simple random walk to zero and it is well known that $Q_\infty(\lambda) = +\infty$ for $\lambda < 0$, while $Q_\infty(\lambda) = 1 - \sqrt{1 - e^{-2\lambda}}$ for $\lambda \geq 0$; see [10]. We point out that $Q_T(\lambda)$ can also be given a closed explicit expression for finite $T$; see Appendix A and, in particular, equation (A.4). Here, it is important to stress that for $T < \infty$, the function $Q_T(\lambda)$ is *analytic and decreasing* on $(\lambda_0^T, +\infty)$, where $\lambda_0^T < 0$ [see equation (A.6)], and $Q_T(\lambda) \to +\infty$ as $\lambda \downarrow \lambda_0^T$, while $Q_T(\lambda) \to 0$ as $\lambda \to \infty$. In particular, when $T < \infty$, the inverse function $(Q_T)^{-1}(\cdot)$ is (analytic and) defined on the whole of $(0, \infty)$, while $(Q_\infty)^{-1}(\cdot)$ is (analytic and) defined only on $(0, 1]$.

THEOREM 1. *Defining $T_\infty = \lim_{N \to \infty} T_N$, the free energy $\phi(\delta, \mathbf{T}) = \phi(\delta, T_\infty)$ depends only on $\delta$ and $T_\infty$, and is given by*

$$(1.10) \qquad \phi(\delta, T_\infty) = \begin{cases} (Q_{T_\infty})^{-1}(e^{-\delta}), & \text{if } T_\infty < +\infty, \\ (Q_\infty)^{-1}(e^{-\delta} \wedge 1), & \text{if } T_\infty = +\infty. \end{cases}$$

*It follows that for $T_\infty < +\infty$, the function $\delta \mapsto \phi(\delta, \mathbf{T})$ is analytic on the whole real line, while for $T_\infty = +\infty$, it is analytic everywhere except at $\delta = 0$.*

So, there are no phase transitions in our model, except in the $T_\infty = +\infty$ case, where $\phi(\delta, \infty)$ is not analytic at $\delta = 0$. This fact is well known because $\phi(\delta, \infty)$ is nothing but the free energy of the classical homogeneous pinning model $\mathbf{P}^\infty_{N,\delta}$; see [12]. In fact, the explicit formula for $Q_\infty(\cdot)$ mentioned above yields

$$(1.11) \qquad \phi(\delta, \infty) = \left(\frac{\delta}{2} - \log\sqrt{2 - e^{-\delta}}\right) \mathbf{1}_{\{\delta \geq 0\}}.$$

Also, in the case where $T_\infty < \infty$, some general properties of $\phi(\delta, T_\infty)$ can be easily derived from Theorem 1, for instance, that $\frac{\partial}{\partial \delta} \phi(\delta, \mathbf{T}) \to 0$ as $\delta \to -\infty$, while $\frac{\partial}{\partial \delta} \phi(\delta, \mathbf{T}) \to \frac{1}{2}$ as $\delta \to +\infty$, both of which have clear physical interpretations, thanks to (1.7).

The proof of Theorem 1 is given in Section 2, using renewal theory ideas. Besides identifying the free energy, we introduce a slightly modified version of the polymer measure $\mathbf{P}^T_{N,\delta}$ which can be given an explicit renewal theory interpretation. This provides a key tool for studying the path behavior (see below).

One consequence of Theorem 1 is that any $\mathbf{T}$ such that $T_\infty = \infty$ yields *the same free energy $\phi(\delta, \mathbf{T}) = \phi(\delta, \infty)$ as the classical homogeneous pinning model*. However, we are going to see that the actual path behavior of $\mathbf{P}^{T_N}_{N,\delta}$ as $N \to \infty$ depends strongly on the speed at which $T_N \to \infty$, a phenomenon which is not captured by the free energy.

1.3. *The scaling behavior.* Henceforth, we focus on the case $\delta > 0$. We assume that $\mathbf{T} = (T_N)_{N \in \mathbb{N}}$ has been chosen such that $T_N \to \infty$ as $N \to \infty$. Then, the free energy $\phi(\delta, \mathbf{T}) = \phi(\delta, \infty)$ is that of the homogeneous pinning model: in particular, $\phi(\delta, \mathbf{T}) > 0$ for every $\delta > 0$. Since $\phi(\delta, \mathbf{T}) = 0$ for $\delta \leq 0$, by convexity and by formula (1.7), it follows that for $\delta > 0$, the typical paths of $\mathbf{P}^{T_N}_{N,\delta}$ touch the interfaces for large $N$ a positive fraction of time and it is customary to say that we are in a *localized regime*.

We now investigate more closely the path properties of $\mathbf{P}^{T_N}_{N,\delta}$. A natural question is: does the polymer visit infinitely many *different* interfaces or does it limit itself to a finite number of them? And, more precisely: what is the scaling behavior of $S_N$ under $\mathbf{P}^{T_N}_{N,\delta}$ as $N \to \infty$?

The answer turns out to depend on the speed at which $T_N \to \infty$. Let $c_\delta$ be the positive constant defined as

$$(1.12) \quad c_\delta := \phi(\delta, \infty) + \log(1 + \sqrt{1 - e^{-2\phi(\delta,\infty)}}) = \frac{\delta}{2} + \log \sqrt{2 - e^{-\delta}},$$

where the right-hand side of (1.12) is obtained with the help of (1.11). Then, the behavior of the sequence $T_N - \frac{1}{c_\delta} \log N$ determines the scaling properties of the polymer measure. More precisely, we have the following result, where $\Longrightarrow$ denotes convergence in law and $\mathcal{N}(0,1)$ the standard normal distribution.

THEOREM 2. *Let $\delta > 0$ and $\mathbf{T} = (T_N)_{N \in \mathbb{N}}$ such that $T_N \to \infty$ as $N \to \infty$.*

(i) *If $T_N - \frac{\log N}{c_\delta} \to -\infty$ as $N \to \infty$, then, under $\mathbf{P}_{N,\delta}^{T_N}$, as $N \to \infty$,*

$$(1.13) \qquad \frac{S_N}{C_\delta(e^{-c_\delta T_N/2} T_N)\sqrt{N}} \Longrightarrow \mathcal{N}(0,1),$$

*where $C_\delta := \sqrt{2e^\delta \phi'(\delta, \infty)} \sqrt{1 - e^{-2\phi(\delta,\infty)}} = (1 - e^{-\delta})\sqrt{\frac{2e^\delta}{2 - e^{-\delta}}}$ is an explicit positive constant.*

(ii) *If there exists $\zeta \in \mathbb{R}$ such that $T_{N'} - \frac{\log N'}{c_\delta} \to \zeta$ along a subsequence $N'$, then, under $\mathbf{P}_{N',\delta}^{T_{N'}}$, as $N' \to \infty$,*

$$(1.14) \qquad \frac{S_{N'}}{T_{N'}} \Longrightarrow S_\Gamma,$$

*where $\Gamma$ is a random variable independent of the $\{S_i\}_{i \geq 0}$ and with a Poisson law of parameter $t_{\delta,\zeta} := 2e^\delta \sqrt{1 - e^{-2\phi(\delta,\infty)}} \phi'(\delta, \infty) \cdot e^{-c_\delta \zeta} = 2e^\delta \frac{(1-e^{-\delta})^2}{2-e^{-\delta}} \cdot e^{-c_\delta \zeta}$.*

(iii) *If $T_N - \frac{\log N}{c_\delta} \to +\infty$ as $N \to \infty$, then the family of laws of $\{S_N\}_{N \in \mathbb{N}}$ under $\mathbf{P}_{N,\delta}^{T_N}$ is tight, that is,*

$$(1.15) \qquad \lim_{L \to \infty} \sup_{N \in \mathbb{N}} \mathbf{P}_{N,\delta}^{T_N}(|S_N| > L) = 0.$$

REMARK 1. It may appear strange that, in point (ii), we have required that $T_{N'} - \frac{\log N'}{c_\delta} \to \zeta$ only along a subsequence $N'$; however, this is simply because $T_N$ takes integer values and therefore the full sequence $T_N - \frac{\log N}{c_\delta}$ cannot have a finite limit. In general, equation (1.14) implies that $S_N/T_N$ is tight when the full sequence $|T_N - \frac{\log N}{c_\delta}|$ is bounded.

The proof of Theorem 2 is distributed across Sections 3, 4 and 5. The crucial idea, described in Section 3.1, is to exploit the renewal theory description given in Section 2. Let us stress the intuitive content of this result. We set $\Delta_N := T_N - \frac{\log(N)}{c_\delta}$ and anticipate that $e^{-c_\delta \Delta_N}$ is the number of jumps made by the polymer of length $N$ between nearest neighbors interfaces. With this in mind, we can provide some more insight into Theorem 2.

- If $\Delta_N \to -\infty$, then the interfaces are departing slowly enough that it is worthwhile for the polymer to visit infinitely many of them. Of course, this is also true when $T_N \equiv T < \infty$ for all $N \in \mathbb{N}$. This situation is not included in Theorem 2 for notational convenience, but a straightforward adaptation of our proof shows that, in this case, $S_N/(\mathcal{C}_T\sqrt{N}) \Longrightarrow \mathcal{N}(0,1)$ for a suitable $\mathcal{C}_T$ satisfying $\mathcal{C}_T \sim C_\delta e^{-c_\delta T/2}T$ as $T \to \infty$, thus according perfectly with (1.13).

    We note that, independently of $(T_N)_N$ (such that $\Delta_N \to -\infty$), the limit law of $S_N$, properly rescaled, is always the standard normal distribution. However, the scaling constants $(e^{-c_\delta T_N/2}T_N)\sqrt{N}$ do depend on the sequence $(T_N)_N$ and, in particular, they are *subdiffusive* as soon as $T_N \to \infty$. Also, note that, by varying $T_N$ from $O(1)$ to the critical case $\frac{\log(N)}{c_\delta} + O(1)$, the scaling constants decrease smoothly from $\sqrt{N}$ to $\log N$.

- If $\Delta_N = O(1)$, then we are in the critical case where the polymer visits a finite number of different interfaces and therefore the scaling behavior of $S_N$ is the same as $T_N$, that is, $S_N \approx \log N$. The explicit form $S_\Gamma$ of the scaling distribution has the following interpretation: the number $\Gamma$ of different interfaces visited by the polymer is distributed according to a Poisson law and, conditionally on $\Gamma$, the polymer just performs $\Gamma$ steps of a simple symmetric random walk on the interfaces.

- If $\Delta_N \to +\infty$, then the only interface visited by the polymer is the $x$-axis. The other interfaces are indeed too distant from the origin to be convenient for the polymer to visit. Therefore, the model $\mathbf{P}_{N,\delta}^{T_N}$ becomes essentially the same as the classical homogeneous pinning model $\mathbf{P}_{N,\delta}^\infty$, where only the interface located at $S = 0$ is present. Since $\delta > 0$, we are in the localized regime for $\mathbf{P}_{N,\delta}^\infty$ and it is well known that $S_N = O(1)$. One could also determine the limit distribution of $S_N$, but we omit this for the sake of conciseness.

As already mentioned, the study of the path behavior in the delocalized regime $\delta < 0$ turns out to be rather different, both from a technical and a physical viewpoint, and will therefore be carried out in a future work.

**2. A renewal theory path to the free energy.** This section is devoted to proving Theorem 1. We also provide a renewal theory description for a slight modification of the polymer measure $\mathbf{P}_{N,\delta}^\infty$, which is the key tool in the following sections.

2.1. *A slight modification.* We consider $\delta \in \mathbb{R}$ and $T \in 2\mathbb{N} \cup \{\infty\}$. It is convenient to introduce the *constrained partition function* $Z_{N,\delta}^{T,c}$, where only the trajectories $(S_i)_i$ that are pinned at an interface at their right extremity are taken into account, that is,

$$(2.1) \qquad Z_{N,\delta}^{T,c} := \mathbf{E}(\exp(H_{N,\delta}^T(S))\mathbf{1}_{\{S_N \in T\mathbb{Z}\}}).$$

In order for the restriction on $\{S_N \in T\mathbb{Z}\}$ to be nontrivial, we work with $Z_{N,\delta}^{T,c}$ only for even $N$. This is the usual parity issue connected with the periodicity of the simple random walk: in fact, $\mathbf{P}(S_N \in T\mathbb{Z}) = 0$ if $N$ is odd (we recall that $T$ was assumed to be even).

The reason for introducing $Z_{N,\delta}^{T,c}$ is that it is easier to handle than the original partition function and, at the same time, it is not too different, as the following lemma shows.

LEMMA 3. *The following relation holds for all $N \in \mathbb{N}, \delta \in \mathbb{R}, T \in 2\mathbb{N}$:*

$$(2.2) \qquad e^{-|\delta|} Z_{2\lfloor N/2 \rfloor,\delta}^{T,c} \leq Z_{N,\delta}^T \leq \sqrt{(N+1) Z_{2N,\delta}^{T,c}}.$$

PROOF. If $N$ is even, then $2\lfloor N/2 \rfloor = N$ and the lower bound in (2.2) follows trivially from the definition (2.1) of $Z_{N,\delta}^{T,c}$. If $N$ is odd, then $2\lfloor N/2 \rfloor = N-1$ and since

$$H_{N,\delta}^T(S) \geq H_{N-1,\delta}^T(S) - |\delta|,$$

the lower bound in (2.2) is proved in full generality.

To prove the upper bound, we observe that by the definition (2.1),

$$Z_{2N,\delta}^{T,c} \geq \mathbf{E}(\exp(H_{2N,\delta}^T(S))\mathbf{1}_{\{S_{2N}=0\}}) = \sum_{k=-N}^{N} \mathbf{E}(\exp(H_{2N,\delta}^T(S))\mathbf{1}_{\{S_N=k\}}\mathbf{1}_{\{S_{2N}=0\}})$$

and, from the Markov property and the time-symmetry $i \mapsto N - i$, we have

$$Z_{2N,\delta}^{T,c} \geq \sum_{k=-N}^{N} [\mathbf{E}(\exp(H_{N,\delta}^T(S))\mathbf{1}_{\{S_N=k\}})]^2.$$

Since $\mathbf{P}(S_N = k) > 0$ if and only if $N$ and $k$ have the same parity, there are only $N+1$ nonzero terms in the sum and, applying Jensen's inequality, we get

$$Z_{2N,\delta}^{T,c} \geq \frac{1}{N+1}\left[\sum_{k=-N}^{N} \mathbf{E}(\exp(H_{N,\delta}^T(S))\mathbf{1}_{\{S_N=k\}})\right]^2 = \frac{1}{N+1}[Z_{N,\delta}^T]^2.$$

Therefore, the upper bound in (2.2) is proved and the proof is complete. □

As a direct consequence of Lemma 3, we observe that to prove the existence of the free energy, that is, of the limit in (1.5), we can safely replace the original partition function $Z_{N,\delta}^{T_N}$ by the constrained one $Z_{N,\delta}^{T_N,c}$, restricting $N$ to the even numbers. The next paragraphs are devoted to obtaining a more explicit expression for $Z_{N,\delta}^{T,c}$.

2.2. *The link with renewal theory.* We start with some definitions. For $T \in 2\mathbb{N} \cup \{\infty\}$, we set $\tau_0^T = 0$ and for $j \in \mathbb{N}$,

$$(2.3) \quad \tau_j^T := \inf\{i \geq \tau_{j-1}^T + 1 : S_i \in T\mathbb{Z}\} \quad \text{and} \quad \varepsilon_j^T := \frac{S_{\tau_j^T} - S_{\tau_{j-1}^T}}{T},$$

where, for $T = \infty$, we agree that $T\mathbb{Z} = \{0\}$. Note that $\tau_j^T$ gives the $j$th epoch at which $S$ touches an interface, while $\varepsilon_j^T$ tells us whether the $j$th interface touched is the same as the $(j-1)$th ($\varepsilon_j^T = 0$), or is the interface above ($\varepsilon_j^T = 1$) or below ($\varepsilon_j^T = -1$). Under the law $\mathbf{P}$ of the simple random walk, we define, for $j = \{0, \pm 1\}$, $n \in \mathbb{N}$ and $\lambda \in \mathbb{R}$, the quantities

$$(2.4) \quad q_T^j(n) := \mathbf{P}(\tau_1^T = n, \varepsilon_1^T = j) \quad \text{and} \quad Q_T^j(\lambda) := \sum_{n=1}^{\infty} e^{-\lambda n} q_T^j(n).$$

Of course, $Q_T^j(\lambda)$ may be (in fact, is) infinite for $\lambda$ negative and large, and, clearly, $q_\infty^{\pm 1}(n) = 0$ for $n \geq 1$ and $Q_\infty^{\pm 1}(\lambda) = 0$ for $\lambda \geq 0$. Note that $q_T^{-1} = q_T^1$ and $Q_T^{-1} = Q_T^1$, so we can focus only on $q_T^j, Q_T^j$ for $j \in \{0, 1\}$. We also set

$$(2.5) \quad \begin{aligned} q_T(n) &:= \sum_{j=0,\pm 1} q_T^j(n) = q_T^0(n) + 2q_T^1(n) = \mathbf{P}(\tau_1^T = n), \\ Q_T(\lambda) &:= \sum_{j=0,\pm 1} Q_T^j(\lambda) = Q_T^0(\lambda) + 2Q_T^1(\lambda) = \mathbf{E}(e^{-\lambda \tau_1^T}). \end{aligned}$$

Next, we introduce

$$(2.6) \qquad \mathcal{H} := \{\mathbb{R} \times 2\mathbb{N}\} \cup \{\mathbb{R}^+ \times \{+\infty\}\}$$

and, for $(\delta, T) \in \mathcal{H}$, we define the quantity $\lambda_{\delta,T}$ by the equation

$$(2.7) \qquad Q_T(\lambda_{\delta,T}) = e^{-\delta}.$$

As will be shown in Appendix A, for $T < \infty$, the function $Q_T(\cdot)$ is analytic and decreasing on $(\lambda_0^T, +\infty)$, with $\lambda_0^T = -\frac{1}{2}\log(1 + (\tan\frac{\pi}{T})^2) < 0$, and such that $Q_T(\lambda) \to +\infty$ as $\lambda \downarrow \lambda_0^T$ and $Q_T(\lambda) \to 0$ as $\lambda \to +\infty$. In particular, equation (2.7) has exactly one solution for every $\delta \in \mathbb{R}$, so $\lambda_{\delta,T}$ is well defined. For $T = \infty$, $Q_T(\cdot)$ is analytic and decreasing on $[0, \infty)$, $Q_T(0) = 1$ and $Q_T(\lambda) \to 0$ as $\lambda \to +\infty$, while $Q_T(\lambda) = \infty$ for $\lambda < 0$. This implies that equation (2.7) has exactly one solution $\lambda_{\delta,\infty}$ for every $\delta \geq 0$ and no solutions for $\delta < 0$. In the next paragraph, we are going to show that when $\lambda_{\delta,T}$ exists, it is nothing but the free energy $\phi(\delta, T)$ (in agreement with Theorem 1).

We are finally ready to introduce, for $(\delta, T) \in \mathcal{H}$, the basic law $\mathcal{P}_{\delta,T}$, under which the sequence of vectors $\{(\xi_i, \varepsilon_i)\}_{i \geq 1}$, taking values in $\mathbb{N} \times \{\pm 1, 0\}$, is i.i.d. with marginal law

$$(2.8) \quad \mathcal{P}_{\delta,T}((\xi_1, \varepsilon_1) = (n, j)) := e^\delta q_T^{|j|}(n) e^{-\lambda_{\delta,T} n}, \qquad n \in \mathbb{N}, j \in \{\pm 1, 0\}.$$

Note that (2.7) ensures that this is indeed a probability law. We then set $\tau_0 = 0$ and $\tau_n = \xi_1 + \cdots + \xi_n$ for $n \geq 1$. We denote by $\tau$ both the sequence of variables $\{\tau_n\}_{n \geq 0}$ and the corresponding random subset of $\mathbb{N} \cup \{0\}$ defined by $\tau = \bigcup_{n \geq 0} \{\tau_n\}$ so that expressions like $\{N \in \tau\}$ make sense. Note that $\{\tau_n\}_{n \geq 0}$ under $\mathcal{P}_{\delta,T}$ is a classical *renewal process* because the increments $\{\tau_n - \tau_{n-1}\}_{n \geq 1} = \{\xi_n\}_{n \geq 1}$ are i.i.d. positive random variables, with law

$$(2.9) \qquad \mathcal{P}_{\delta,T}(\tau_1 = n) := e^\delta q_T(n) e^{-\lambda_{\delta,T} n}, \qquad n \in \mathbb{N}.$$

Because of the periodicity of the simple random walk, $q_T(n) = 0$ for all odd $n \in \mathbb{N}$ and $q_T(n) > 0$ for all even $n \in \mathbb{N}$ (we recall that we are only considering the case of even $T$). Therefore, the renewal process is *periodic with period* 2.

We now have all of the ingredients needed to give an explicit expression for the partition function in terms of the jumps made by $S$ between interfaces. This can be done for $(\delta, T) \in \mathcal{H}$ and for $Z_{N,\delta}^{T,c}$ [recall (2.1)], as follows. For $k, n \in \mathbb{N}$, $k \leq n$, we define the set

$$\mathcal{S}_{k,n} := \{t \in (\mathbb{N} \cup \{0\})^{k+1} : 0 = t_0 < t_1 < \cdots < t_k = n\}.$$

Then, for $\lambda \in \mathbb{R}$ and $N$ even, we can write

$$(2.10) \quad \begin{aligned} Z_{N,\delta}^{T,c} &= \sum_{k=1}^{N} \sum_{\sigma \in \{-1,0,1\}^k} \sum_{t \in \mathcal{S}_{k,N}} \prod_{l=1}^{k} e^\delta q_T^{|\sigma_l|}(t_l - t_{l-1}) \\ &= e^{\lambda N} \sum_{k=1}^{N} \sum_{\sigma \in \{-1,0,1\}^k} \sum_{t \in \mathcal{S}_{k,N}} \prod_{l=1}^{k} e^\delta q_T^{|\sigma_l|}(t_l - t_{l-1}) e^{-\lambda(t_l - t_{l-1})}. \end{aligned}$$

Then, setting $\lambda = \lambda_{\delta,T}$ and recalling (2.8), we can rewrite (2.10) as

$$(2.11) \qquad Z_{N,\delta}^{T,c} = e^{\lambda_{\delta,T} \cdot N} \mathcal{P}_{\delta,T}(N \in \tau).$$

We stress that this equation retains a crucial importance in our approach. In fact, the behavior of $Z_{N,\delta}^{T,c}$ is reduced to the asymptotic properties of the renewal process $\tau$.

The next step is to lift relation (2.11) from the constrained partition function to the *constrained polymer measure* $\mathbf{P}_{N,\delta}^{T,c}$, defined for even $N$ as

$$\mathbf{P}_{N,\delta}^{T,c}(\cdot) := \mathbf{P}_{N,\delta}^{T}(\cdot | S_N \in T\mathbb{Z}).$$

Recalling the definition (1.6) of $L_{N,T}$, for $(\delta, T) \in \mathcal{H}$, $k \leq N$, $t \in \mathcal{S}_{k,N}$ and $\sigma \in \{\pm 1, 0\}^k$, in analogy to (2.10), we can write

$$
\begin{aligned}
(2.12) \quad & \mathbf{P}_{N,\delta}^{T,c}(L_{N,T} = k, (\tau_i^T, \varepsilon_i^T) = (t_i, \sigma_i), 1 \leq i \leq k) \\
& = \frac{e^{\lambda_{\delta,T} N}}{Z_{N,\delta}^{T,c}} \prod_{l=1}^{k} e^{\delta} q_T^{|\sigma_l|}(t_l - t_{l-1}) e^{-\lambda_{\delta,T}(t_l - t_{l-1})}.
\end{aligned}
$$

Therefore, from (2.8) and (2.11), we obtain

$$
\begin{aligned}
(2.13) \quad & \mathbf{P}_{N,\delta}^{T,c}(L_{N,T} = k, (\tau_i^T, \varepsilon_i^T) = (t_i, \sigma_i), 1 \leq i \leq k) \\
& = \mathcal{P}_{\delta,T}(L_N = k, (\tau_i, \varepsilon_i) = (t_i, \sigma_i), 1 \leq i \leq k | N \in \tau),
\end{aligned}
$$

where $L_N := \sup\{j \geq 1 : \tau_j \leq N\}$, in analogy with (1.6). Thus, the process $\{(\tau_i^T, \varepsilon_i^T)\}_{i > L_N}$ under $\mathbf{P}_{N,\delta}^{T,c}$ is distributed like $\{(\tau_i, \varepsilon_i)\}_{i > L_N}$ under the explicit law $\mathcal{P}_{\delta,T}$, conditioned on the event $\{N \in \tau\}$. The crucial point is that $\{\tau_i\}_i$ under $\mathcal{P}_{\delta,T}$ is a genuine *renewal process*. This fact is the key to the path results that we prove in the next section because we will show that the constrained law $\mathbf{P}_{N,\delta}^{T,c}$ is not too different from the original law $\mathbf{P}_{N,\delta}^{T}$.

2.3. *Proof of Theorem 1.* Thanks to Lemma 3, to prove Theorem 1, it suffices to show that for every sequence $(T_N)_N$ such that $T_N \to T_\infty$ as $N \to \infty$, we have

$$
(2.14) \quad \lim_{N \to \infty, N \text{ even}} \frac{1}{N} \log Z_{\delta,N}^{T_N,c} = \begin{cases} (Q_{T_\infty})^{-1}(e^{-\delta}), & \text{if } T_\infty < \infty, \\ (Q_{T_\infty})^{-1}(e^{-\delta} \wedge 1), & \text{if } T_\infty = \infty, \end{cases}
$$

where we recall that $Q_T(\cdot)$ was introduced in (1.9). Recall, also, that for $(\delta, T) \in \mathcal{H}$, we have $(Q_T)^{-1}(e^{-\delta}) = \lambda_{\delta,T}$ [see (2.7)].

First, consider the case where $T_\infty < \infty$, that is, $T_\infty \in \mathbb{N}$. The sequence $(T_N)_N$ then eventually takes the constant value $T_N = T_\infty$ and, thanks to (2.11) and (2.7), we can write

$$
(2.15) \quad \frac{1}{N} \log Z_{\delta,N}^{T_\infty,c} = (Q_{T_\infty})^{-1}(e^{-\delta}) + \frac{1}{N} \log \mathcal{P}_{\delta,T_\infty}(N \in \tau).
$$

It therefore remains to show that the last term in the right-hand side vanishes as $N \to \infty$, $N$ even, and we are done [as a byproduct, we also show that $\lambda_{\delta,T_\infty}$ coincides with the free energy $\phi(\delta, T_\infty)$]. We recall that the process $\tau = \{\tau_n\}_n$ under $\mathcal{P}_{\delta,T_\infty}$ is a classical renewal process with step-mean

$$
(2.16) \quad m(\delta, T_\infty) := \mathcal{E}_{\delta,T_\infty}(\tau_1) < +\infty.
$$

The fact that $m(\delta, T_\infty) < +\infty$ is easily checked by (2.9) because, by construction, $\lambda_{\delta,T_\infty} > \lambda_0^{T_\infty}$; see (2.5), (2.7) and the following lines. Since the renewal process $(\{\tau_n\}_n, \mathcal{P}_{\delta,T_\infty})$ has period 2, the renewal theorem yields

$$
(2.17) \quad \lim_{N \to \infty, N \text{ even}} \mathcal{P}_{\delta,T_\infty}(N \in \tau) = \frac{2}{m(\delta, T_\infty)} > 0
$$

and, looking back to (2.15), we see that (2.14) is proved.

Next, we consider the case where $T_\infty = +\infty$, that is, $T_N \to +\infty$ as $N \to \infty$. We can rewrite equation (2.15) as

$$(2.18) \qquad \frac{1}{N} \log Z_{\delta,N}^{T_N,c} = (Q_{T_N})^{-1}(e^{-\delta}) + \frac{1}{N} \log \mathcal{P}_{\delta,T_N}(N \in \tau).$$

We start by considering the first term in the right-hand side of (2.18), by proving the following lemma.

LEMMA 4. *For every $\delta \in \mathbb{R}$,*

$$(2.19) \qquad \lim_{T \to \infty, T \in 2\mathbb{N}} (Q_T)^{-1}(e^{-\delta}) = (Q_\infty)^{-1}(e^{-\delta} \wedge 1).$$

PROOF. To this end, we observe that as $T \to \infty$, the variable $\tau_1^T$, defined in (2.3), converges a.s. toward $\tau_1^\infty := \inf\{i > 0 : S_i = 0\}$, that is, the first return to zero of the simple random walk. Accordingly, by dominated convergence (or by direct verification), $Q_T(\lambda)$ converges as $T \to \infty$, *for every* $\lambda \in [0, +\infty)$, toward $Q_\infty(\lambda) = 1 - \sqrt{1 - e^{-2\lambda}}$. Since $Q_\infty(\cdot)$ is strictly decreasing, it is easily checked that the inverse functions also converge, that is, for every $y \in (0, 1]$, we have $(Q_T)^{-1}(y) \to (Q_\infty)^{-1}(y)$ as $T \to \infty$ so that (2.19) is checked for $\delta \geq 0$. On the other hand, when $\delta < 0$, we have $\lambda_0^T < (Q_T)^{-1}(e^{-\delta}) < 0$ because, as we already mentioned, $Q_T(\cdot)$ is decreasing and $Q_T(\lambda) \to \infty$ as $\lambda \downarrow \lambda_0^T$ and $Q_T(0) = 1$. Moreover, $\lambda_0^T$ vanishes as $T \to \infty$ [see (A.6)] and, consequently, $(Q_T)^{-1}(e^{-\delta}) \to 0$ as $T \to \infty$. Hence, (2.19) also holds for $\delta < 0$. $\square$

Using Lemma 4 and the fact that $\mathcal{P}_{\delta,T_N}(N \in \tau) \leq 1$, by (2.18), we obtain

$$\limsup_{N \to \infty, N \text{ even}} \frac{1}{N} \log Z_{\delta,N}^{T_N,c} \leq (Q_\infty)^{-1}(e^{-\delta} \wedge 1).$$

Hence, to complete the proof of (2.14), it remains to show that for every $\delta \in \mathbb{R}$,

$$(2.20) \qquad \liminf_{N \to \infty, N \text{ even}} \frac{1}{N} \log Z_{\delta,N}^{T_N,c} \geq (Q_\infty)^{-1}(e^{-\delta} \wedge 1).$$

We start by considering the case where $\delta \leq 0$, hence $(Q_\infty)^{-1}(e^{-\delta} \wedge 1) = 0$. We give a very rough lower bound on $Z_{\delta,N}^{T_N,c}$, namely, for $N$ even, we can write

$$(2.21) \qquad \begin{aligned} Z_{\delta,N}^{T_N,c} &\geq \mathbf{E}(\exp(H_{N,\delta}^{T_N}(S))\mathbf{1}_{\{S_i \notin T_N \mathbb{Z}, \forall 1 \leq i \leq N-1\}} \mathbf{1}_{\{S_N = 0\}}) \\ &= e^\delta \cdot q_{T_N}^0(N), \end{aligned}$$

where we recall that $q^0_{T_N}(N) = \mathbf{P}(\tau_1^{T_N} = N; S_N = 0)$ was defined in (2.4). (If $N$ is odd, then the same formula holds, just replacing $N$ by $N-1$, and the following considerations are easily adapted.) So, we are left with showing that $q^0_{T_N}(N)$ does not decay exponentially fast as $N \to \infty$: by the explicit formula (A.7), we have

$$q^0_{T_N}(N) \geq \frac{2}{T_N} \cos^{N-2}\left(\frac{\pi}{T_N}\right) \sin^2\left(\frac{\pi}{T_N}\right).$$

At this stage, by using the fact that $\sin^2(x) \sim x^2$ as $x \to 0$, we can assert that for $N$ large enough, $\sin(\pi/T_N) \geq \pi/(2T_N)$ and since, by assumption, $T_N \leq N$, we obtain

$$\mathbf{P}_1(\tau_1^{T_N} = N-1; S_{N-1} = 0) \geq \frac{\pi^2}{2N^3} e^{(N-2)\log \cos(\pi/T_N)},$$

which, by (2.21), shows that (2.20) holds [note that the right-hand side of (2.20) is zero for $\delta \leq 0$].

Finally, we have to prove that equation (2.20) holds true for $\delta > 0$. By (2.18) and Lemma 4, it suffices to show that

(2.22) $$\liminf_{N \to \infty} \frac{1}{N} \log \mathcal{P}_{\delta, T_N}(N \in \tau) = 0.$$

This is not straightforward because the law $\mathcal{P}_{\delta, T_N}$ changes with $N$ and therefore some uniformity is needed. Let us be more precise: by the renewal theorem [see (2.17)], *for fixed $T$*, we have that, as $n \to \infty$ along the even numbers,

$$\mathcal{P}_{\delta, T}(n \in \tau) \longrightarrow \frac{2}{m(\delta, T)},$$

where $m(\delta, T)$ was introduced in (2.16). At the same time, as $T \to \infty$, we have

$$m(\delta, T) \longrightarrow m(\delta, \infty),$$

as will be proven in Lemma 6 below. Since $T_N \to \infty$ as $N \to \infty$, the last two equations suggest that for $N$ large, $\mathcal{P}_{\delta, T_N}(N \in \tau)$ should be close to $2/m(\delta, \infty)$. To show that this is indeed the case, we are going to apply Theorem 2 in [14], which is a uniform version of the renewal theorem. First, recall that, by Lemma 4, $\lambda_{\delta, T} \to \lambda_{\delta, \infty} > 0$ as $T \to \infty$, $T \in 2\mathbb{N}$ and, moreover, $\lambda_{\delta, T} > 0$ for every $T \in 2\mathbb{N}$. Hence, there exist $C_1, C_2 > 0$ such that $C_1 \leq \lambda_{\delta, T} \leq C_2$ for every $T \in 2\mathbb{N}$. We are ready to verify the following two conditions:

(1) when $\delta > 0$ is fixed and $T$ varies in $2\mathbb{N}$, the family of renewal processes $(\{\tau_n\}_n, \mathcal{P}_{\delta, T})$ restricted to the even numbers is uniformly aperiodic, in the sense of Definition 1 in [14], because $\mathcal{P}_{\delta, T}(\tau_1 = 2) = e^\delta q_T(2) e^{-2\lambda_{\delta, T}} \geq (e^\delta/2) \cdot e^{-2C_2} > 0$ for all $T \in 2\mathbb{N}$;

(2) when $\delta > 0$ is fixed and $T$ varies in $2\mathbb{N}$, the family of renewal processes $(\{\tau_n\}_n, \mathcal{P}_{\delta,T})$ has uniformly summable tails, in the sense of Definition 2 in [14], because

$$\mathcal{P}_{\delta,T}(\tau_1 \geq t) \leq \sum_{r=t}^{\infty} e^{-C_1 r} = \frac{e^{-C_1 t}}{1 - e^{-C_1}}.$$

We can therefore apply Theorem 2 in [14], which yields the following lemma. This implies (2.22) and the proof of Theorem 1 is thus complete.

LEMMA 5. *Fix $\delta > 0$. Then, for every $\varepsilon > 0$, there exist $N_0 \in \mathbb{N}$ such that for every $T \in 2\mathbb{N}$ and for all $N \geq N_0$, $N$ even, we have*

$$\left| \mathcal{P}_{\delta,T}(N \in \tau) - \frac{2}{m(\delta, \infty)} \right| \leq \varepsilon.$$

LEMMA 6. *For all $\delta > 0$ and $k \in \mathbb{N}$,*

(2.23) $$\lim_{T \to \infty} \mathcal{E}_{\delta,T}((\tau_1)^k) = \mathcal{E}_{\delta,\infty}((\tau_1)^k).$$

PROOF. By Lemma 4, we know that for $\delta > 0$, we have $\lambda_{\delta,T} \to \lambda_{\delta,\infty} > 0$ as $T \to \infty$, $T \in 2\mathbb{N}$. Thus, by writing

$$\mathcal{E}_{\delta,T}((\tau_1)^k) = e^\delta \sum_{n=1}^{\infty} n^k q_T(n) e^{-\lambda_{\delta,T} n},$$

it suffices to apply the dominated convergence theorem [since $q_T(n) \leq 1$]. □

REMARK 2. Now that we have proven that the free energy $\phi(\delta, T)$ indeed equals the right-hand side of (2.14), we can restate Lemma 4 in the following way:

(2.24) $$\lim_{T \to \infty} \phi(\delta, T) = \phi(\delta, \infty) \qquad \forall \delta \in \mathbb{R}.$$

REMARK 3. For $(\delta, T) \in \mathcal{H}$, we know that $\lambda_{\delta,T} = \phi(\delta, T)$. Consequently, we will use $\phi(\delta, T)$ instead of $\lambda_{\delta,T}$ in what follows.

**3. Proof of Theorem 2(i).** This section is devoted to the proof of part (i) of Theorem 2. We recall that $\delta > 0$ is fixed and that $T_N - \frac{1}{c_\delta} \log N \to -\infty$ as $N \to \infty$, where $c_\delta$ is defined in (1.12).

We recall that $(\tau_i^T, \varepsilon_i^T)_{i \geq 1}$ defined in (2.3) under $\mathbf{P}_{N,\delta}^{T_N}$ represents the jump process of the polymer between the interfaces, whereas $(\tau_i, \varepsilon_i)_{i \geq 1}$ introduced in (2.8) under the law $\mathcal{P}_{\delta,T_N}$ represents an auxiliary renewal process. For $N \geq 1$, we set

$$Y_N^T = \sum_{i=1}^{N} \varepsilon_i^T \quad \text{and recall that}$$

(3.1)
$$L_{N,T} = \sup\{j \geq 1 : \tau_j^T \leq N\}.$$

Analogously, we set

(3.2) $$Y_N = \sum_{i=1}^{N} \varepsilon_i \quad \text{and recall that} \quad L_N = \sup\{j \geq 1 : \tau_j \leq N\}.$$

3.1. *General strategy.* Let us describe the strategy of our proof. The aim is to determine the asymptotic behavior of $S_N$ under $\mathbf{P}_{N,\delta}^{T_N}$ as $N \to \infty$. The starting point is given by the following considerations:

- by definition, we have $S_N = T \cdot Y_{L_{N,T}}^T + O(T)$, hence the behavior of $S_N$ can be recovered from that of $L_{N,T}$ and $\{Y_n^T\}_n$;
- it turns out that the *free* polymer measure $\mathbf{P}_{N,\delta}^{T_N}$ is not too different from the *constrained* one $\mathbf{P}_{N,\delta}^{T_N,c} = \mathbf{P}_{N,\delta}^{T_N}(\cdot | S_N \in T_N \mathbb{Z})$, which, in turn, is closely linked to the law $\mathcal{P}_{\delta,T_N}$ introduced in Section 2.2; see, in particular, (2.13).

For these reasons, the first part of the proof of Theorem 2 consists of determining the asymptotic behavior of $\{Y_n\}_n$ and $L_N$ under $\mathcal{P}_{\delta,T_N}$. This is carried out in Section 3.3 (Step 1) and Section 3.4 (Step 2) below, exploiting ideas and techniques from random walks and renewal theory. The second part of the proof is devoted to showing that the law $\mathcal{P}_{\delta,T_N}$ can indeed be replaced by $\mathbf{P}_{N,\delta}^{T_N,c}$ [see Section 3.5 (Step 3)] and, finally, by $\mathbf{P}_{N,\delta}^{T_N}$; see Section 3.6 (Step 4).

Let us give a closer (heuristic) look at the core of the proof. *For fixed $T$,* the process $\{Y_n\}_n$ under $\mathcal{P}_{\delta,T}$ is just a symmetric random walk on $\mathbb{Z}$ with step law

$$\mathcal{P}_{\delta,T}(Y_1 = j) = \mathcal{P}_{\delta,T}(\varepsilon_1 = j) = e^\delta Q_T^{|j|}(\phi(\delta,T)), \qquad j \in \{\pm 1, 0\};$$

see equations (2.8) and (2.4), (2.5). In particular, the central limit theorem yields

(3.3) $$Y_N \approx C_T \sqrt{N} \qquad \text{under } \mathcal{P}_{\delta,T} \text{ as } N \to \infty,$$

where $C_T = \sqrt{2e^\delta Q_T^1(\phi(\delta,T))}$ is the standard deviation of $Y_1$.

Of course, we are interested in the case where $T = T_N$ is no longer fixed, but varies with $N$, more precisely $T_N \to \infty$ as $N \to \infty$. It is then easy to see that $C_{T_N} \to 0$. However, if it happens that $C_{T_N} \sqrt{N} \to \infty$ as $N \to \infty$, one may hope that equation (3.3) still holds with $T$ replaced by $T_N$. This is indeed true, as we are going to show. To determine the asymptotic behavior of $C_T$, the following lemma is useful.

LEMMA 7. *Fix $\delta > 0$. Then, as $T \to \infty$,*

$$Q_T^1(\phi(\delta,T)) = \sqrt{1 - e^{-2\phi(\delta,\infty)}} e^{-c_\delta T}(1 + o(1)), \tag{3.4}$$

*where $c_\delta = \phi(\delta,\infty) + \log(1 + \sqrt{1 - e^{-2\phi(\delta,\infty)}})$ [recall (1.12)].*

This shows that the condition $C_{T_N}\sqrt{N} \to \infty$ as $N \to \infty$ is equivalent to $T_N - \frac{\log N}{c_\delta} \to -\infty$, which is exactly the hypothesis of part (i) of Theorem 2. As we mentioned, in this case, we show that (3.3) still holds, so

$$Y_N \approx C_{T_N}\sqrt{N} \approx C^* e^{-c_\delta T_N/2}\sqrt{N} \qquad \text{under } \mathcal{P}_{\delta,T_N} \text{ as } N \to \infty, \tag{3.5}$$

with $C^* = \sqrt{2e^\delta \sqrt{1 - e^{-2\phi(\delta,\infty)}}}$.

Now, let us return to $S_N$. By definition, we have $S_N = T_N \cdot Y_{L_{N,T_N}}^{T_N} + O(T_N)$ and, from equation (1.7), we get $L_{N,T_N} \approx cN$, with $c = \phi'(\delta,\infty) > 0$. Moreover, as was already mentioned, the law $\mathcal{P}_{\delta,T_N}$ can be replaced by the original polymer measure $\mathbf{P}_{N,\delta}^{T_N}$ without changing the asymptotic behavior. Together with (3.5), these considerations yield

$$S_N \approx T_N \cdot Y_{cN}^{T_N} \approx C_\delta(e^{-c_\delta T_N/2} T_N)\sqrt{N} \qquad \text{under } \mathbf{P}_{N,\delta}^{T_N} \text{ as } N \to \infty,$$

where $C_\delta := C^*\sqrt{c} = \sqrt{2e^\delta \phi'(\delta,\infty)\sqrt{1 - e^{-2\phi(\delta,\infty)}}}$. Note that this accords exactly with the result of Theorem 2.

PROOF OF LEMMA 7. We can rewrite the second relation in (A.3) as

$$\begin{aligned}Q_T^1(\lambda) &= \sqrt{1 - e^{-2\lambda}} \cdot e^{-\widetilde{c}_\lambda T} \\ &\quad \cdot \frac{1}{1 + ((1 - \sqrt{1 - e^{-2\lambda}})/(1 + \sqrt{1 - e^{-2\lambda}}))^T},\end{aligned} \tag{3.6}$$

where $\widetilde{c}_\lambda := \lambda + \log(1 + \sqrt{1 - e^{-2\lambda}})$. We have to replace $\lambda$ by $\phi(\delta,T)$ in this relation and study the asymptotic behavior as $T \to \infty$.

Observe that $\phi(\delta,T)$ and $\phi(\delta,\infty)$ are both strictly positive since $\delta > 0$ and, moreover, $\phi(\delta,T) \to \phi(\delta,\infty)$ as $T \to \infty$ (see Remark 2). This easily implies that the last factor in the right-hand side of (3.6) is $1 + o(1)$, hence, as $T \to \infty$,

$$Q_T^1(\phi(\delta,T)) = \sqrt{1 - e^{-2\phi(\delta,\infty)}} e^{-\widetilde{c}_{\phi(\delta,T)}T}(1 + o(1)).$$

To prove (3.4), it remains to show that $\widetilde{c}_{\phi(\delta,T)}T = c_\delta T + o(1)$ as $T \to \infty$. Since $c_\delta = \widetilde{c}_{\phi(\delta,\infty)}$, this follows once we show that $|\phi(\delta,T) - \phi(\delta,\infty)| = o(\frac{1}{T})$.

To this end, we fix $\varepsilon > 0$ such that $\phi(\delta,T) \geq \varepsilon$ for every $T$. By equation (A.4), there exists $\kappa = \kappa_\varepsilon > 0$ such that, *uniformly for* $\lambda \in [\varepsilon,\infty)$,

$$Q_T(\lambda) = 1 - \sqrt{1 - e^{-2\lambda}} + O(e^{-\kappa T}) \qquad (T \to \infty).$$

Recalling that $Q_\infty(\lambda) = 1 - \sqrt{1 - e^{-2\lambda}}$ and that, by Theorem 1, $e^{-\delta} = Q_T(\phi(\delta,T)) = Q_\infty(\phi(\delta,\infty))$, we obtain

$$Q_\infty(\phi(\delta,T)) - Q_\infty(\phi(\delta,\infty)) = O(e^{-\kappa T}) \qquad (T \to \infty).$$

Since $Q_\infty(\lambda)$ is continuously differentiable with nonzero derivative for $\lambda > 0$, it follows that $\phi(\delta,T) - \phi(\delta,\infty) = O(e^{-\kappa T})$ and the proof is complete. □

3.2. *Preparation.* We start the proof of Theorem 2 by rephrasing equation (1.13), which is our goal, in a slightly different form. We recall that $T_N - \frac{1}{c_\delta} \log N \to -\infty$ as $N \to \infty$ or, equivalently, $e^{-c_\delta T_N} N \to \infty$, and that, by construction, $|S_N - Y^{T_N}_{L_{N,T_N}} \cdot T_N| \leq T_N$. Therefore, equation (1.13) is equivalent to the following: for all $x \in \mathbb{R}$,

$$(3.7) \qquad \lim_{N \to \infty} \mathbf{P}^{T_N}_{N,\delta}\left(\frac{Y^{T_N}_{L_{N,T_N}}}{C_\delta \sqrt{e^{-c_\delta T_N} N}} \leq x\right) = P(\mathcal{N}(0,1) \leq x),$$

where

$$(3.8) \qquad C_\delta = \sqrt{2e^\delta \phi'(\delta,\infty)} \sqrt{1 - e^{-2\phi(\delta,\infty)}}.$$

Recall the definition (2.9) of the renewal process $(\tau, \mathcal{P}_{\delta,T})$. For $\delta > 0$ and $T \in 2\mathbb{N} \cup \{+\infty\}$, we set

$$(3.9) \qquad s_T := \frac{1}{\mathcal{E}_{\delta,T}(\tau_1)} \in (0,\infty).$$

Differentiating the relation $Q_T(\phi(\delta,T)) = e^{-\delta}$, one obtains $\phi'(\delta,T) = -e^{-\delta}/Q'_T(\phi(\delta,T))$ and, by direct computation,

$$(3.10) \quad \begin{aligned} \mathcal{E}_{\delta,T}(\tau_1) &= e^\delta \sum_{n \in \mathbb{N}} n q_T(n) e^{-\phi(\delta,T)n} \\ &= -e^\delta Q'_T(\phi(\delta,T)) = \frac{1}{\phi'(\delta,T)} \qquad \forall T \in 2\mathbb{N}. \end{aligned}$$

In particular, $\phi'(\delta,\infty) = s_\infty$. Recalling Lemma 7 and setting $Q^1_{T_N} := Q^1_{T_N}(\phi(\delta,T_N))$ for conciseness, we can finally restate (3.7) as

$$(3.11) \qquad \lim_{N \to \infty} \mathbf{P}^{T_N}_{N,\delta}\left(\frac{Y^{T_N}_{L_{N,T_N}}}{\sqrt{s_\infty}\sqrt{2e^\delta Q^1_{T_N} N}} \leq x\right) = P(\mathcal{N}(0,1) \leq x),$$

which is exactly what we are going to prove. This will be achieved in four steps. We stress that the assumption $T_N - \frac{1}{c_\delta} \log N \to -\infty$ as $N \to \infty$ is equivalent to $Q^1_{T_N} \cdot N \to \infty$.

3.3. *Step 1.* In this step, we consider the auxiliary renewal process of law $\mathcal{P}_{\delta,T_N}$ and prove that, for $x \in \mathbb{R}$,

$$(3.12) \qquad \lim_{N \to \infty} \mathcal{P}_{\delta,T_N}\left(\frac{Y_N}{\sqrt{2e^\delta Q^1_{T_N} N}} \leq x\right) = P(\mathcal{N}(0,1) \leq x).$$

Under the law $\mathcal{P}_{\delta,T_N}$, $(\varepsilon_1, \ldots, \varepsilon_N)$ are symmetric i.i.d. random variables taking values $-1, 0, 1$. Therefore, they satisfy

$$(3.13) \qquad \mathcal{E}_{\delta,T_N}(|\varepsilon_1|^3) = \mathcal{E}_{\delta,T_N}((\varepsilon_1)^2) = 2e^\delta Q^1_{T_N}$$

and we can apply the Berry–Esseen theorem, which gives

$$(3.14) \qquad \left|\mathcal{P}_{\delta,T_N}\left(\frac{Y_N}{\alpha_\delta(N,T_N)} \leq x\right) - P(\mathcal{N}(0,1) \leq x)\right| \leq \frac{3\mathcal{E}_{\delta,T_N}(|\varepsilon_1|^3)}{\mathcal{E}_{\delta,T_N}(\varepsilon_1^2)^{3/2}\sqrt{N}}$$
$$= \frac{3}{\sqrt{2e^\delta Q^1_{T_N} N}}.$$

Since $Q^1_{T_N} \cdot N \to \infty$ by assumption, equation (3.12) is proved.

3.4. *Step 2.* In this step, we prove that, for $x \in \mathbb{R}$,

$$(3.15) \qquad \lim_{N \to \infty} \mathcal{P}_{\delta,T_N}\left(\frac{Y_{L_N}}{\sqrt{s_\infty}\sqrt{2e^\delta Q^1_{T_N} N}} \leq x\right) = P(\mathcal{N}(0,1) \leq x).$$

The idea is to show that $L_N \approx s_\infty \cdot N$ and then to apply (3.12). We need the following lemma.

LEMMA 8. *For every $\varepsilon > 0$, there exists $T_0 = T_0(\varepsilon) \in \mathbb{N}$ such that*

$$(3.16) \qquad \lim_{N \to \infty} \sup_{T \geq T_0} \mathcal{P}_{\delta,T}\left(\left|\frac{L_N}{N} - s_\infty\right| > \varepsilon\right) = 0.$$

PROOF. Lemma 6 yields $s_T \to s_\infty$ as $T \to \infty$ [we recall the definition (3.9)]. Therefore, we fix $T_0 = T_0(\varepsilon)$ such that $|s_\infty - s_T| \leq \frac{\varepsilon}{2}$ for $T \geq T_0$ and consequently

$$\mathcal{P}_{\delta,T}\left(\left|\frac{L_N}{N} - s_\infty\right| > \varepsilon\right) \leq \mathcal{P}_{\delta,T}\left(\left|\frac{L_N}{N} - s_T\right| > \frac{\varepsilon}{2}\right).$$

Setting $\xi_i = \tau_i - \tau_{i-1}$ and $\widetilde{\xi}_i = \xi_i - \frac{1}{s_T}$, by Chebyshev's inequality, we get

$$\mathcal{P}_{\delta,T}\left(\frac{L_N}{N} > s_T + \varepsilon\right) = \mathcal{P}_{\delta,T}(\tau_{\lfloor (s_T+\varepsilon)N \rfloor} \leq N)$$

$$= \mathcal{P}_{\delta,T}\left(-\widetilde{\xi}_1 - \cdots - \widetilde{\xi}_{\lfloor (s_T+\varepsilon)N \rfloor} \geq \frac{\varepsilon N}{s_T}\right)$$

$$\leq \frac{s_T^2(s_T + \varepsilon)\mathcal{E}_{\delta,T}(\widetilde{\xi}_1^2)}{\varepsilon^2 N}.$$

By Lemma 6, both the sequences $T \mapsto s_T$ and $T \mapsto \mathcal{E}_{\delta,T}(\widetilde{\xi}_1^2)$ are bounded and therefore the right-hand side above vanishes as $N \to \infty$, uniformly in $T$. The event $\{\frac{L_N}{N} < s_T - \varepsilon\}$ is treated with analogous arguments and the proof is then complete.  □

We set

$$\frac{Y_{L_N}}{\sqrt{s_\infty}\sqrt{2e^\delta Q_{T_N}^1 N}} = \frac{Y_{\lfloor s_\infty N \rfloor}}{\sqrt{s_\infty}\sqrt{2e^\delta Q_{T_N}^1 N}} + \frac{Y_{L_N} - Y_{\lfloor s_\infty N \rfloor}}{\sqrt{s_\infty}\sqrt{2e^\delta Q_{T_N}^1 N}} =: V_N + G_N.$$

Step 1 [see equation (3.12)] implies directly that $V_N$ converges in law to $\mathcal{N}(0,1)$. Therefore, it remains to prove that $G_N$ converges in probability to 0. For $\eta, \varepsilon > 0$, we write

$$\mathcal{P}_{\delta,T_N}(|G_N| > \eta)$$

(3.17) $$\leq \mathcal{P}_{\delta,T_N}\left(|G_N| > \eta, \left|\frac{L_N}{N} - s_\infty\right| \leq \varepsilon\right) + \mathcal{P}_{\delta,T_N}\left(\left|\frac{L_N}{N} - s_\infty\right| > \varepsilon\right)$$

$$\leq \mathcal{P}_{\delta,T_N}(U_{\varepsilon,N} > \eta) + \mathcal{P}_{\delta,T_N}\left(\left|\frac{L_N}{N} - s_\infty\right| > \varepsilon\right).$$

3.5. *Step 3.* This is the most delicate step, where we show that one can replace the free measure $\mathcal{P}_{\delta,T_N}$ by the constrained one, $\mathcal{P}_{\delta,T_N}(\cdot | N \in \tau)$. More precisely, we prove that, for $x \in \mathbb{R}$,

(3.18)
$$\lim_{N \to \infty, N \text{ even}} \mathcal{P}_{\delta,T_N}\left(\frac{Y_{L_N}}{\sqrt{s_\infty}\sqrt{2e^\delta Q_{T_N}^1 N}} \leq x \Big| N \in \tau\right)$$
$$= P(\mathcal{N}(0,1) \leq x).$$

We note that one can safely replace $L_N$ with $L_{N-\lfloor \sqrt{T_N} \rfloor}$ in the left-hand side because $Y_{L_{N-\lfloor \sqrt{T_N} \rfloor}}$ differs from $Y_{L_N}$ by at most $\pm 1$. The same is true for equation (3.15), which we rewrite for convenience as follows:

(3.19) $$\lim_{N \to \infty} \mathcal{P}_{\delta,T_N}\left(\frac{Y_{L_{N-\lfloor \sqrt{T_N} \rfloor}}}{\sqrt{s_\infty}\sqrt{2e^\delta Q_{T_N}^1 N}} \leq x\right) = P(\mathcal{N}(0,1) \leq x).$$

By summing over the locations of the last point $t$ in $\tau$ before $N - \lfloor\sqrt{T_N}\rfloor$ and of the first point $r$ in $\tau$ after $N - \lfloor\sqrt{T_N}\rfloor$, and using the Markov property, we obtain

$$\mathcal{P}_{\delta,T_N}\left(\frac{Y_{L_{N-\lfloor\sqrt{T_N}\rfloor}}}{\sqrt{s_\infty}\sqrt{2e^\delta Q^1_{T_N}N}} \leq x \Big| N \in \tau\right)$$

$$= \frac{1}{\mathcal{P}_{\delta,T_N}(N \in \tau)}$$

$$\times \sum_{t=0}^{N-\lfloor\sqrt{T_N}\rfloor} \sum_{r=t+1}^{t+\lfloor\sqrt{T_N}\rfloor} \mathcal{P}_{\delta,T_N}\left(\frac{Y_{L_{N-\lfloor\sqrt{T_N}\rfloor}}}{\sqrt{s_\infty}\sqrt{2e^\delta Q^1_{T_N}N}} \leq x, N - \lfloor\sqrt{T_N}\rfloor - t \in \tau\right)$$

$$\cdot \mathcal{P}_{\delta,T_N}(\tau_1 = r) \cdot \mathcal{P}_{\delta,T_N}(t + \lfloor\sqrt{T_N}\rfloor - r \in \tau).$$

Introducing the function

$$\Theta_{\delta,N}(t) := \frac{\sum_{r=t+1}^{t+\lfloor\sqrt{T_N}\rfloor} \mathcal{P}_{\delta,T_N}(\tau_1 = r) \cdot \mathcal{P}_{\delta,T_N}(t + \lfloor\sqrt{T_N}\rfloor - r \in \tau)}{\mathcal{P}_{\delta,T_N}(N \in \tau) \cdot \sum_{r=t+1}^{\infty} \mathcal{P}_{\delta,T_N}(\tau_1 = r)},$$

we can write

(3.20)
$$\mathcal{P}_{\delta,T_N}\left(\frac{Y_{L_{N-\lfloor\sqrt{T_N}\rfloor}}}{\sqrt{s_\infty}\sqrt{2e^\delta Q^1_{T_N}N}} \leq x \Big| N \in \tau\right)$$
$$= \sum_{t=0}^{N-\lfloor\sqrt{T_N}\rfloor} \mathcal{P}_{\delta,T_N}\left(\frac{Y_{L_{N-\lfloor\sqrt{T_N}\rfloor}}}{\sqrt{s_\infty}\sqrt{2e^\delta Q^1_{T_N}N}} \leq x, N - \lfloor\sqrt{T_N}\rfloor - t \in \tau\right)$$
$$\cdot \mathcal{P}_{\delta,T_N}(\tau_1 > t) \cdot \Theta_{\delta,N}(t).$$

Note that if we set $\Theta_{\delta,N}(t) \equiv 1$, the right-hand side of the last relation becomes the left-hand side of (3.19). In fact, $\Theta_{\delta,N}(t)$ is nothing but the Radon–Nikodym derivative of the conditioned law $\mathcal{P}_{\delta,T_N}(\cdot|N \in \tau)$ with respect to the free one $\mathcal{P}_{\delta,T_N}$. We are going to show that $\Theta_{\delta,N}(t) \to 1$ as $N \to \infty$, *uniformly in the values of $t$ that have the same parity as* $\lfloor\sqrt{T_N}\rfloor$ [otherwise, $\Theta_{\delta,N}(t) = 0$]. If we succeed in this, equation (3.18) will follow from (3.19).

Let us set $K_N(n) := \mathcal{P}_{\delta,T_N}(\tau_1 = n)$ and $u_N(n) := \mathcal{P}_{\delta,T_N}(n \in \tau)$ so that we can rewrite $\Theta_{\delta,N}(t)$ as

(3.21) $$\Theta_{\delta,N}(t) := \frac{\sum_{r=t+1}^{t+\lfloor\sqrt{T_N}\rfloor} K_N(r) \cdot u_N(t + \lfloor\sqrt{T_N}\rfloor - r)}{u_N(N) \cdot \sum_{r=t+1}^{\infty} K_N(r)}.$$

We recall that

$$K_N(n) = e^\delta e^{-\phi(\delta,T_N)\cdot n} q_{T_N}(n)$$

[see (2.9)] and $q_T(\cdot)$ is defined in (2.4). We are going to show the following: *for every $\varepsilon > 0$, there exists $N_0 = N_0(\varepsilon)$ such that for every $N \geq N_0$ and for all values of $t \leq N - \lfloor\sqrt{T_N}\rfloor$ that have the same parity as $\lfloor\sqrt{T_N}\rfloor$, we have*

$$(3.22) \qquad 1 - \varepsilon \leq \Theta_{\delta,N}(t) \leq 1 + \varepsilon.$$

The proof of this step will then be complete. First, we need a preliminary lemma.

LEMMA 9. *For every $\eta > 0$, there exists $N_1 = N_1(\eta)$ such that for every $N \geq N_1$ and for all $0 \leq t \leq N - \lfloor\sqrt{T_N}\rfloor$, we have*

$$(3.23) \qquad \sum_{r=t+\lfloor\sqrt{T_N}\rfloor/2}^{\infty} K_N(r) \leq \eta \cdot \left(\sum_{r=t+1}^{t+\lfloor\sqrt{T_N}\rfloor/2} K_N(r)\right).$$

PROOF. First, we observe that, by the explicit formulae in (A.7), the following upper bound holds for every $T, n \in \mathbb{N}$ with $n \geq 2$:

$$\max\{q_T^0(n), 2q_T^1(n)\} \leq \frac{2}{T} \sum_{\nu=1}^{\lfloor(T-1)/2\rfloor} \cos^{n-2}\left(\frac{\pi\nu}{T}\right) \sin^2\left(\frac{\pi\nu}{T}\right).$$

We can bound the left-hand side of (3.23) as

$$\sum_{r=t+\lfloor\sqrt{T_N}\rfloor/2}^{\infty} K_N(r) \leq e^{\delta} e^{-\phi(\delta,T_N)\cdot(t+\lfloor\sqrt{T_N}\rfloor/2)} \sum_{r=t+\lfloor\sqrt{T_N}\rfloor/2}^{\infty} q_{T_N}(r)$$

and, since $q_T(r) = q_T^0(r) + 2q_T^1(r)$, we have

$$\sum_{r=t+\lfloor\sqrt{T_N}\rfloor/2}^{\infty} q_{T_N}(r) \leq 2 \sum_{r=t+\lfloor\sqrt{T_N}\rfloor/2}^{\infty} \left(\frac{2}{T_N} \sum_{\nu=1}^{\lfloor(T_N-1)/2\rfloor} \cos^{r-2}\left(\frac{\pi\nu}{T_N}\right) \sin^2\left(\frac{\pi\nu}{T_N}\right)\right)$$

$$= \frac{4}{T_N} \sum_{\nu=1}^{\lfloor(T_N-1)/2\rfloor} \frac{(\cos(\pi\nu/T_N))^{t-2+\lfloor\sqrt{T_N}\rfloor/2}}{1-\cos(\pi\nu/T_N)} \sin^2\left(\frac{\pi\nu}{T_N}\right)$$

$$\leq \frac{4}{T_N} \cdot \frac{T_N}{2} \left(\cos\left(\frac{\pi}{T_N}\right)\right)^{t-2+\lfloor\sqrt{T_N}\rfloor/2} \cdot 2,$$

where we have used the fact that $\sin^2 x/(1-\cos x) = 1 + \cos x \leq 2$ for $x \in (0, \frac{\pi}{2}]$. Therefore,

$$\sum_{r=t+\lfloor\sqrt{T_N}\rfloor/2}^{\infty} K_N(r) \leq 4 e^{\delta} e^{-\phi(\delta,T_N)\cdot(t+\lfloor\sqrt{T_N}\rfloor/2)} \left(\cos\left(\frac{\pi}{T_N}\right)\right)^{t-2+\lfloor\sqrt{T_N}\rfloor/2}.$$
(3.24)

Next, we bound from below the right-hand side of (3.23):

$$\sum_{r=t+1}^{t+\lfloor\sqrt{T_N}\rfloor/2} K_N(r) \geq e^\delta e^{-\phi(\delta,T_N)\cdot(t+2)}(q_{T_N}^0(t+1) + q_{T_N}^0(t+2)).$$

One of the two numbers $t+1, t+2$ is even—call it $\ell$. We can then apply equation (A.7) to get

$$q_{T_N}^0(\ell) = \frac{2}{T_N} \sum_{\nu=1}^{\lfloor(T_N-1)/2\rfloor} \cos^{\ell-2}\left(\frac{\pi\nu}{T_N}\right) \sin^2\left(\frac{\pi\nu}{T_N}\right)$$

$$\geq \frac{2}{T_N} \cos^{\ell-2}\left(\frac{\pi}{T_N}\right) \sin^2\left(\frac{\pi}{T_N}\right),$$

hence

$$(3.25) \quad \sum_{r=t+1}^{t+\lfloor\sqrt{T_N}\rfloor/2} K_N(r) \geq e^\delta e^{-\phi(\delta,T_N)\cdot(t+2)} \frac{2}{T_N} \cos^t\left(\frac{\pi}{T_N}\right) \sin^2\left(\frac{\pi}{T_N}\right).$$

The ratio of the right-hand sides of equations (3.24) and (3.25) equals

$$2T_N e^{-\phi(\delta,T_N)\cdot(\lfloor\sqrt{T_N}\rfloor/2-2)} \frac{(\cos(\pi/T_N))^{\lfloor\sqrt{T_N}\rfloor/2-2}}{\sin^2(\pi/T_N)}$$

$$\leq \frac{8}{\pi^2}(T_N)^3 e^{-\phi(\delta,T_N)\cdot(\lfloor\sqrt{T_N}\rfloor/2-2)}.$$

Since the right-hand side no longer depends on $t$ and vanishes as $N \to \infty$, the proof is complete. □

Let us return to the proof of (3.22). We first observe that, thanks to Lemma 5, for every $\eta > 0$, there exists $N_2 = N_2(\eta)$ and that for all $N \in \mathbb{N}$ and all $r \geq N_2$, $r$ even, we have

$$(1-\eta)2s_\infty \leq u_N(r) \leq (1+\eta)2s_\infty$$

[$s_\infty$ is defined in (3.9)]. Henceforth, we assume that $t$ has the same parity as $\lfloor\sqrt{T_N}\rfloor$. Then, if $N$ is large, such that $\lfloor\sqrt{T_N}\rfloor/2 \geq N_2$, we can bound $\Theta_{\delta,N}(t)$ [recall (3.21)] by

$$\Theta_{\delta,N}(t) \leq \frac{(1+\eta)2s_\infty \sum_{r=t+1}^{t+\lfloor\sqrt{T_N}\rfloor/2} K_N(r) + \sum_{t+\lfloor\sqrt{T_N}\rfloor/2+1}^{t+\lfloor\sqrt{T_N}\rfloor} K_N(r)}{(1-\eta)2s_\infty \sum_{r=t+1}^{t+\lfloor\sqrt{T_N}\rfloor/2} K_N(r)}$$

and, if $N \geq N_1$, we can apply Lemma 9 to obtain

$$\Theta_{\delta,N}(t) \leq \frac{1+\eta+\eta/(2s_\infty)}{1-\eta} \leq 1+\varepsilon,$$

provided that $\eta$ is chosen sufficiently small. Therefore, the upper bound in (3.22) is proved. The lower bound is analogous: for large $N$, we have

$$\Theta_{\delta,N}(t) \geq \frac{(1-\eta)2s_\infty \sum_{r=t+1}^{t+\lfloor\sqrt{T_N}\rfloor/2} K_N(r)}{(1+\eta)2s_\infty \sum_{r=t+1}^{t+\lfloor\sqrt{T_N}\rfloor/2} K_N(r) + \sum_{t+\lfloor\sqrt{T_N}\rfloor/2+1}^{t+\lfloor\sqrt{T_N}\rfloor} K_N(r)}$$

and, again applying Lemma 9, we finally obtain

$$\Theta_{\delta,N}(t) \geq \frac{1-\eta}{1+\eta+\eta/(2s_\infty)} \geq 1-\varepsilon,$$

provided $\eta$ is small. Recalling (3.20) and the following lines, the step is completed.

3.6. *Step 4.* In this step, we finally complete the proof of Theorem 2(i), proving equation (3.11), that we rewrite for convenience: for every $x \in \mathbb{R}$,

$$(3.26) \qquad \lim_{N\to\infty} \mathbf{P}_{N,\delta}^{T_N}\left(\frac{Y_{L_{N,T_N}}^{T_N}}{\sqrt{s_\infty}\sqrt{2e^\delta Q_{T_N}^1 N}} \leq x\right) = P(\mathcal{N}(0,1) \leq x).$$

We start by summing over the location $\mu_N := \tau_{L_{N,T_N}}^{T_N}$ of the last point in $\tau^{T_N}$ before $N$ (we henceforth assume that $N$ is even):

$$\mathbf{P}_{N,\delta}^{T_N}\left(\frac{Y_{L_{N,T_N}}^{T_N}}{\sqrt{s_\infty}\sqrt{2e^\delta Q_{T_N}^1 N}} \leq x\right) = \sum_{\ell=0}^{N} \mathbf{P}_{N,\delta}^{T_N}\left(\frac{Y_{L_{N,T_N}}^{T_N}}{\sqrt{s_\infty}\sqrt{2e^\delta Q_{T_N}^1 N}} \leq x \Big| \mu_N = N-\ell\right)$$

$$\cdot \mathbf{P}_{N,\delta}^{T_N}(\mu_N = N-\ell).$$

Of course, only the terms with $\ell$ even are nonzero. We start by showing that we can truncate the sum at a finite number of terms. To this end, we estimate

$$\mathbf{P}_{N,\delta}^{T_N}(\mu_N = N-\ell) = \frac{\mathbf{E}(\exp(H_{N-\ell,\delta}^{T_N}(S))\mathbf{1}_{\{N-\ell\in\tau\}}) \cdot \mathbf{P}(\tau_1 > \ell)}{\mathbf{E}(\exp(H_{N,\delta}^{T_N}(S)))}.$$

We focus on the denominator: inserting the event $\{N-\ell \in \tau\}$ and using the Markov property yields

$$\mathbf{E}(\exp(H_{N,\delta}^{T_N}(S))) \geq \mathbf{E}(\exp(H_{N-\ell,\delta}^{T_N}(S))\mathbf{1}_{\{N-\ell\in\tau\}}) \cdot \mathbf{E}(\exp(H_{\ell,\delta}^{T_N}(S))),$$

hence

$$\mathbf{P}_{N,\delta}^{T_N}(\mu_N = N-\ell) \leq \frac{\mathbf{P}(\tau_1 > \ell)}{\mathbf{E}(\exp(H_{\ell,\delta}^{T_N}(S)))} \leq \frac{1}{\mathbf{E}(\exp(H_{\ell,\delta}^\infty(S)))} = \frac{1}{Z_{\ell,\delta}^\infty},$$

where we have used the elementary fact that $\mathbf{E}(\exp(H_{\ell,\delta}^T(S))) \geq \mathbf{E}(\exp(H_{\ell,\delta}^\infty(S)))$ for every $T \in \mathbb{N}$; see (1.1) and (1.3). Note that the right-hand side above no longer depends on $N$ and that $Z_{\ell,\delta}^\infty \asymp \exp(\phi(\delta,\infty) \cdot \ell)$ as $\ell \to \infty$, where $\asymp$ denotes equivalence in the Laplace sense; see [12]. Since $\phi(\delta,\infty) > 0$ for $\delta > 0$, it follows that for every $\varepsilon > 0$, there exists $\ell_0 = \ell_0(\varepsilon)$ such that for every $N \in \mathbb{N}$, we have

$$(3.27) \qquad \sum_{\ell=\ell_0+1}^{N} \mathbf{P}_{N,\delta}^{T_N}(\mu_N = N - \ell) \leq \varepsilon.$$

As a consequence, we have

$$\left| \mathbf{P}_{N,\delta}^{T_N}\left( \frac{Y_{L_{N,T_N}}^{T_N}}{\sqrt{s_\infty}\sqrt{2e^\delta Q_{T_N}^1 N}} \leq x \right) \right.$$
$$\left. - \sum_{\ell=0}^{\ell_0} \mathbf{P}_{N,\delta}^{T_N}\left( \frac{Y_{L_{N,T_N}}^{T_N}}{\sqrt{s_\infty}\sqrt{2e^\delta Q_{T_N}^1 N}} \leq x \middle| \mu_N = N - \ell \right) \cdot \mathbf{P}_{N,\delta}^{T_N}(\mu_N = N - \ell) \right|$$
$$\leq \varepsilon.$$

Therefore, to complete the proof of (3.26), it remains to show that for every fixed $\ell \in \mathbb{N} \cup \{0\}$,

$$(3.28) \quad \lim_{N \to \infty} \mathbf{P}_{N,\delta}^{T_N}\left( \frac{Y_{L_{N,T_N}}^{T_N}}{\sqrt{s_\infty}\sqrt{2e^\delta Q_{T_N}^1 N}} \leq x \middle| \mu_N = N - \ell \right) = P(\mathcal{N}(0,1) \leq x).$$

However, this is easy. In fact, on the event $\{\mu_N = N - \ell\}$, we have $Y_{L_{N,T_N}}^{T_N} = Y_{L_{N-\ell,T_N}}^{T_N}$ and, by the Markov property, we get

$$\mathbf{P}_{N,\delta}^{T_N}\left( \frac{Y_{L_{N,T_N}}^{T_N}}{\sqrt{s_\infty}\sqrt{2e^\delta Q_{T_N}^1 N}} \leq x \middle| \mu_N = N - \ell \right)$$
$$= \mathbf{P}_{N,\delta}^{T_N}\left( \frac{Y_{L_{N-\ell,T_N}}^{T_N}}{\sqrt{s_\infty}\sqrt{2e^\delta Q_{T_N}^1 N}} \leq x \middle| N - \ell \in \tau \right).$$

However, arguing as in Section 2.2 [see, in particular, (2.13)], we have that

$$\mathbf{P}_{N,\delta}^{T_N}\left( \frac{Y_{L_{N-\ell,T_N}}^{T_N}}{\sqrt{s_\infty}\sqrt{2e^\delta Q_{T_N}^1 N}} \leq x \middle| N - \ell \in \tau \right)$$
$$= \mathcal{P}_{\delta,T_N}\left( \frac{Y_{L_{N-\ell}}}{\sqrt{s_\infty}\sqrt{2e^\delta Q_{T_N}^1 N}} \leq x \middle| N - \ell \in \tau \right).$$

Therefore, (3.28) follows easily from (3.18).

**4. Proof of Theorem 2(ii).** This section is devoted to the proof of part (ii) of Theorem 2, which, in a sense, is the *critical regime*. We stress that $\delta > 0$ is fixed throughout the section. The assumption in part (ii) is that the sequence $(T_N)_N$ is such that $T_{N'} - \frac{\log N'}{c_\delta} \to \zeta$ along a subsequence $N'$, where $\zeta \in \mathbb{R}$ (the reason for considering only a subsequence is explained in Remark 1). However, for notational convenience, in this section, we drop the subsequence and assume that for some $\zeta \in \mathbb{R}$, as $N \to \infty$,

(4.1)
$$T_N - \frac{\log N}{c_\delta} \longrightarrow \zeta \quad \text{or, equivalently,}$$
$$Q^1_{T_N} \cdot N \longrightarrow \sqrt{1 - e^{-2\phi(\delta,\infty)}} e^{-c_\delta \zeta},$$

where we have used Lemma 7 and we recall the shorthand $Q^1_{T_N} := Q^1_{T_N}(\phi(\delta, T_N))$ introduced in the previous section.

We recall that the variables $(\xi_i, \varepsilon_i, \tau_i)_{i \geq 1}$ are defined under the law $\mathcal{P}_{\delta, T_N}$ [see (2.8)]. We now introduce the successive epochs $(\theta_i)_{i \geq 0}$ at which the jump process changes interface, by setting $\theta_0 = 0$ and, for $j \geq 1$,

(4.2) $\quad \theta_j := \inf\{m > \theta_{j-1} : \exists i \in \mathbb{N} \text{ such that } \tau_i = m \text{ and } |\varepsilon_i| = 1\}.$

The number of these jumps occurring before time $N$ is given by

(4.3) $\quad L'_N := \sup\{j \geq 0 : \theta_j \leq N\} = \#\{i \leq L_N : |\varepsilon_i| = 1\}.$

Note that $\theta \subseteq \tau$, where, as usual, we identify $\theta = \{\theta_n\}_n$ with a (random) subset of $\mathbb{N} \cup \{0\}$.

We split the proof into three steps.

4.1. *Step 1.* We start by proving that under $\mathcal{P}_{\delta, T_N}$, the variable $L'_N$ converges in law to a Poisson law of parameter $t_{\delta,\zeta}$ with $t_{\delta,\zeta} := 2e^\delta \sqrt{1 - e^{-2\phi(\delta,\infty)}} \phi'(\delta,\infty) \cdot e^{-c_\delta \zeta}$, that is,

(4.4) $\quad \lim_{N \to \infty} \mathcal{P}_{\delta, T_N}(L'_N = j) = e^{-t_{\delta,\zeta}} \frac{(t_{\delta,\zeta})^j}{j!} \quad \forall j \in \mathbb{N} \cup \{0\}.$

We note that $\{|\varepsilon_i|\}_{i \geq 1}$ under $\mathcal{P}_{\delta, T_N}$ is a sequence of i.i.d. Bernoulli trials with success probability given by

(4.5) $\quad p_{T_N} := \mathcal{P}_{\delta, T_N}(|\varepsilon_1| = 1) = 2e^\delta Q^1_{T_N}.$

We also set

$$\Delta := \inf\{i \geq 1 : |\varepsilon_i| = 1\}.$$

Note that $(\theta_j - \theta_{j-1})_{j \geq 1}$ are i.i.d. random variables. Moreover, we can write

$$\theta_1 = \sum_{j=1}^{\Delta} \xi_j.$$

We now study the asymptotic behavior of $\theta_j$ and, by (4.3), we derive that of $L'_N$. The building blocks are given in the following lemma.

LEMMA 10. *The following convergences in law hold as $N \to \infty$, under $\mathcal{P}_{\delta,T_N}$:*

(4.6)
$$\frac{\xi_\Delta}{N} \Longrightarrow 0, \qquad \frac{\Delta - 1}{N} \Longrightarrow \mathrm{Exp}(v_\delta),$$
$$\frac{1}{\Delta - 1} \sum_{j=1}^{\Delta-1} \xi_j \Longrightarrow \mathcal{E}_{\delta,\infty}(\xi_1),$$

*where $v_{\delta,\zeta} := 2e^\delta \sqrt{1 - e^{-2\phi(\delta,\infty)}} e^{-c_\delta \zeta}$ and $\mathrm{Exp}(\lambda)$ denotes the exponential law of parameter $\lambda$, that is, $P(\mathrm{Exp}(\lambda) \in \mathrm{d}x) = \lambda e^{-\lambda x} \mathbf{1}_{\{x \geq 0\}} \mathrm{d}x$.*

PROOF. For the first relation, it suffices to show that $\mathcal{E}_{\delta,T_N}(\xi_\Delta/N)$ vanishes as $N \to \infty$. By definition, the variable $\xi_\Delta$ gives the length of a jump conditioned to occur between two different interfaces, namely, $\xi_\Delta$ has the same law as $\xi_1$, conditionally on the event $\{|\varepsilon_1| = 1\}$. This leads to the following formula [see (2.8)]:

(4.7)
$$\mathcal{E}_{\delta,T_N}\left(\frac{\xi_\Delta}{N}\right) = \frac{1}{Q^1_{T_N} N} \sum_{n=1}^{\infty} n q^1_{T_N}(n) e^{-\phi(\delta,T_N)n}.$$

By (4.1), $Q^1_{T_N} N \to c' > 0$ as $N \to \infty$ and, for every fixed $n \geq 1$, we observe that, plainly, $q^1_{T_N}(n) \to 0$ as $N \to \infty$ [in fact, $q^1_T(n) = 0$ for $T > n$]. Since $\phi(\delta,T_N) \to \phi(\delta,\infty) > 0$ as $N \to \infty$ (see Remark 2), by dominated convergence, the right-hand side of (4.7) vanishes as $N \to \infty$.

For the second relation in (4.6), note that the variable $\Delta$ has a geometric law of parameter $p_{T_N}$, that is, for all $j \in \mathbb{N}$,

$$\mathcal{P}_{\delta,T_N}(\Delta = j) = (1 - p_{T_N})^{j-1} p_{T_N}.$$

Since $N \cdot p_{T_N} \to v_{\delta,\zeta} = 2e^\delta \sqrt{1 - e^{-2\phi(\delta,\infty)}} e^{-c_\delta \zeta}$ as $N \to \infty$ [see (4.1) and (4.5)], it is well known (and easy to check) that $\Delta/N$ converges to an exponential law of parameter $v_\delta$ and, of course, the same is true for $(\Delta - 1)/N$.

Next, we focus on the third relation in (4.6). Since $\mathcal{P}_{\delta,T_N}(\Delta \leq \sqrt{N}) \to 0$ as $N \to \infty$, by the result just proven, it suffices to consider, for $\varepsilon > 0$, the quantity

(4.8)
$$\mathcal{P}_{\delta,T_N}\left(\left|\frac{1}{\Delta-1}\sum_{j=1}^{\Delta-1}\xi_j - \mathcal{E}_{\delta,\infty}(\xi_1)\right| > \varepsilon, \Delta > \sqrt{N}\right)$$
$$= \sum_{l=\lceil\sqrt{N}\rceil}^{\infty} \mathcal{P}_{\delta,T_N}(\Delta = l) \mathcal{P}_{\delta,T_N}\left(\left|\frac{1}{l-1}\sum_{j=1}^{l-1}\xi_j - \mathcal{E}_{\delta,\infty}(\xi_1)\right| > \varepsilon \Big| \Delta = l\right).$$

To evaluate the last term, we note that under $\mathcal{P}_{\delta,T_N}(\cdot|\Delta = l)$, the variables $\xi_1, \ldots, \xi_{l-1}$ are i.i.d. with marginal law simply given by the law of $\xi_1$, conditionally on the event $\{\varepsilon_1 = 0\}$ (which means that the jump occurs at the same interface). Denoting this law, for simplicity, by $\mathcal{P}^0_{\delta,T_N}$, we have, for $n \geq 1$,

$$(4.9) \qquad \mathcal{P}^0_{\delta,T_N}(\xi_1 = n) = \frac{1}{1 - 2e^\delta Q^1_{T_N}} q^0_{T_N}(n) e^\delta e^{-\phi(\delta, T_N)n}.$$

By (4.1), we have $Q^1_{T_N} \to 0$ as $n \to \infty$. Moreover, $q^0_{T_N}(n) \to q_\infty(n)$, by definition, and $\phi(\delta, T_N) \to \phi(\delta, \infty) > 0$, by Remark 2. These considerations yield, by dominated convergence, $\mathcal{E}^0_{\delta,T_N}(\xi_1) \to \mathcal{E}_{\delta,\infty}(\xi_1)$ and $\text{Var}^0_{\delta,T_N}(\xi_1) \to \text{Var}_{\delta,\infty}(\xi_1)$ as $N \to \infty$. In particular, in the right-hand side of (4.8), we can replace $\mathcal{E}_{\delta,\infty}(\xi_1)$ by $\mathcal{E}^0_{\delta,T_N}(\xi_1)$ and $\varepsilon$ by (say) $\varepsilon/2$, and we get an upper bound for large $N$. Applying Chebyshev's inequality, we obtain

$$(4.10) \quad \mathcal{P}_{\delta,T_N}\left(\left|\frac{1}{l-1}\sum_{j=1}^{l-1} \xi_j - \mathcal{E}^0_{\delta,T_N}(\xi_1)\right| > \frac{\varepsilon}{2}\bigg|\Delta = l\right) \leq \frac{4\,\text{Var}^0_{\delta,T_N}(\xi_1)}{\varepsilon^2(l-1)}.$$

This shows that the right-hand side of (4.8) vanishes as $N \to \infty$ and this completes the proof. □

By writing

$$\frac{\theta_1}{N} = \frac{\Delta - 1}{N} \cdot \frac{1}{\Delta - 1}\sum_{j=1}^{\Delta-1} \xi_j + \frac{\xi_\Delta}{N}$$

and applying Lemma 10, we can easily conclude that $\theta_1/N$ converges in law to an exponential distribution of parameter $t_{\delta,\zeta}$ given by

$$t_{\delta,\zeta} := v_{\delta,\zeta}/\mathcal{E}_{\delta,\infty}(\xi_1) = 2e^\delta\sqrt{1 - e^{-2\phi(\delta,\infty)}}e^{-c_\delta\zeta} \cdot \phi'(\delta,\infty),$$

having used (3.10). By independence, for every fixed $j \in \mathbb{N}$, the variable $\theta_j/N$ converges to a gamma law with parameters $(j, t_{\delta,\zeta})$, hence, by (4.3), the variable $L'_N$ converges to a Poisson law of parameter $t_{\delta,\zeta}$. This completes the step.

4.2. *Step 2.* In this step, we want to prove that under the law $\mathcal{P}_{\delta,T_N}(\cdot|N \in \tau)$, with $N \in 2\mathbb{N}$, the quantity $L'_N$ still converges to a Poisson distribution of parameter $t_{\delta,\zeta}$, that is,

$$(4.11) \quad \lim_{N \to \infty, N \text{ even}} \mathcal{P}_{\delta,T_N}(L'_N = j|N \in \tau) = e^{-t_{\delta,\zeta}}\frac{(t_{\delta,\zeta})^j}{j!} \qquad \forall j \in \mathbb{N} \cup \{0\}.$$

We start by elaborating a little on (4.4). Fix $L \in \mathbb{N}$ and write, by the renewal property,

$$\mathcal{P}_{\delta,T_N}(\tau \cap (N-L,N] = \varnothing) = \sum_{r=0}^{N-L} \sum_{s=N+1}^{\infty} u_N(r) \cdot K_N(s-r),$$

where we recall the definitions $u_N(n) := \mathcal{P}_{\delta,T_N}(n \in \tau)$ and $K_N(n) := \mathcal{P}_{\delta,T_N}(\tau_1 = n)$. Since $u_N(r) \leq 1$ and $K_N(n) \leq e^{\delta} e^{-\phi(\delta,T_N)\cdot n}$ [see (2.9)], and since $\phi(\delta, T_N) \to \phi(\delta, \infty) > 0$ as $N \to \infty$ (see Remark 2), it follows that

$$\mathcal{P}_{\delta,T_N}(\tau \cap (N-L,N] = \varnothing) \leq e^{\delta} \sum_{r=0}^{N-L} \sum_{s=N+1}^{\infty} e^{-\phi(\delta,T_N)\cdot(s-r)}$$

(4.12)
$$\leq C \cdot e^{-C' \cdot L},$$

where $C, C'$ are suitable positive constants depending only on $\delta$. This means that the probability of the event $\{\tau \cap (N-L,N] = \varnothing\}$ can be made arbitrarily small, *uniformly in $N$*, by taking $L$ large. It is then easy to see that equation (4.4) yields the following: for all $\varepsilon > 0$ and $j \in \mathbb{N} \cup \{0\}$, there exist $N_0, L_0$ such that for all $N \geq N_0$ and $L \geq L_0$, we have

(4.13)
$$\mathcal{P}_{\delta,T_N}(L'_N = j | \tau \cap (N-L,N] \neq \varnothing)$$
$$\in \left(e^{-t_{\delta,\zeta}} \frac{(t_{\delta,\zeta})^j}{j!} - \varepsilon, e^{-t_{\delta,\zeta}} \frac{(t_{\delta,\zeta})^j}{j!} + \varepsilon\right).$$

Next, we show that equation (4.11) follows from (4.13). The idea is that conditioning on the event $\{\tau \cap (N-L,N] \neq \varnothing\}$, that is, that there is a renewal epoch in $(N-L,N]$, is the same as conditioning on $\{N - i \in \tau\}$ for some $i = 0, \ldots, L-1$ and the latter is essentially independent of $i$. More precisely, we have the following lemma.

LEMMA 11. *For every $i \in 2\mathbb{N} \cup \{0\}$, the following relation holds as $N \to \infty$, with $N \in 2\mathbb{N}$:*

(4.14) $\quad \mathcal{P}_{\delta,T_N}(L'_N = j | N \in \tau) = \mathcal{P}_{\delta,T_N}(L'_{N-i} = j | N - i \in \tau) + \varepsilon_i(N),$

*where $\varepsilon_i(N) \to 0$ as $N \to \infty$.*

PROOF. Note that $\{L'_N = j\} = \{\theta_j \leq N, \theta_{j+1} > N\}$. First, we restrict the expectation on the event $\{\theta_j \leq N - \sqrt{N}\}$, which has almost full probability. In fact, for fixed $i \in 2\mathbb{N} \cup \{0\}$,

(4.15)
$$\mathcal{P}_{\delta,T_N}(L'_{N-i} = j, \theta_j > N - \sqrt{N} | N - i \in \tau)$$
$$\leq \frac{\mathcal{P}_{\delta,T_N}(N - \sqrt{N} < \theta_j \leq N)}{\mathcal{P}_{\delta,T_N}(N - i \in \tau)} = o(1)$$

as $N \to \infty$, $N \in 2\mathbb{N}$, because $\theta_j/N$ converges as $N \to \infty$ to an atom-free law (in fact, to a gamma law) by Step 1 and, by Lemma 5, $\mathcal{P}_{\delta,T_N}(N-i \in \tau) \to 2/m(\delta,\infty) > 0$ as $N \to \infty$. Specializing (4.15) to $i=0$, we can therefore write as $N \to \infty$, $N \in 2\mathbb{N}$,

$$\mathcal{P}_{\delta,T_N}(L'_N = j | N \in \tau) = \mathcal{P}_{\delta,T_N}(L'_N = j, \theta_j \le N - \sqrt{N} | N \in \tau) + o(1)$$
$$= \mathcal{P}_{\delta,T_N}(\theta_j \le N - \sqrt{N}, \theta_{j+1} > N | N \in \tau) + o(1).$$

The renewal property then yields

(4.16)
$$\mathcal{P}_{\delta,T_N}(L'_N = j | N \in \tau)$$
$$= \sum_{r=1}^{\lfloor N-\sqrt{N} \rfloor} \mathcal{P}_{\delta,T_N}(\theta_j = r) \cdot \frac{\mathcal{P}_{\delta,T_N}(\theta_1 > N-r, N-r \in \tau)}{\mathcal{P}_{\delta,T_N}(N \in \tau)} + o(1).$$

We now study the term $\mathcal{P}_{\delta,T_N}(\theta_1 > l, l \in \tau)$. We have

$$\mathcal{P}_{\delta,T_N}(\theta_1 > l, l \in \tau)$$
(4.17)
$$= \sum_{k=1}^{l} \sum_{0=:t_0 < t_1 < \cdots < t_k = l} \prod_{j=1}^{k} e^{\delta} q^0_{T_N}(t_j - t_{j-1}) e^{-\phi(\delta,T_N)(t_j - t_{j-1})}$$
$$= e^{\nu_N \cdot l} \cdot \sum_{0=:t_0 < t_1 < \cdots < t_k = l} \prod_{j=1}^{k} \widetilde{K}^0_N(t_j - t_{j-1}),$$

where we have set, for $n \in \mathbb{N}$,

$$\widetilde{K}^0_N(n) := e^{\delta} q^0_{T_N}(n) e^{-(\phi(\delta,T_N)+\nu_N) \cdot n}$$

and we fix $\nu_N < 0$ such that $\sum_{n \in \mathbb{N}} \widetilde{K}^0_N(n) = 1$, that is, $Q^0_{T_N}(\phi(\delta,T_N)+\nu_N) = e^{-\delta}$, which is always possible because $Q^0_T(\lambda)$ diverges as $\lambda \downarrow \lambda^0_T$; see Appendix A. Denoting by $\widetilde{\mathcal{P}}^0_{\delta,T_N}$ the global law of $\tau$, when the step distribution is $\widetilde{K}^0_N(n)$, we can rewrite (4.17) with $l = N-r$ as

(4.18) $\mathcal{P}_{\delta,T_N}(\theta_1 > N-r, N-r \in \tau) = e^{\nu_N \cdot (N-r)} \cdot \widetilde{\mathcal{P}}^0_{\delta,T_N}(N-r \in \tau).$

Plainly, as $N \to \infty$, we have $q^0_{T_N}(n) \to q_\infty(n)$ for every $n \in \mathbb{N}$, where we recall that $q_\infty(n)$ is the return time distribution for the simple random walk; see Section 2.2. Hence, $\nu_N \to 0$ and $\widetilde{K}^0_N(n) \to \mathcal{P}_{\delta,\infty}(n \in \tau)$ as $N \to \infty$. Then, a slight modification of Lemma 5 shows that, for any fixed $r \in 2\mathbb{N}$, $\widetilde{\mathcal{P}}^0_{\delta,T_N}(N-r \in \tau) \to 2/m(\delta,\infty) > 0$ as $N \to \infty$. Then, in equation (4.18), we can replace $N$ by $N-i$, any fixed $i \in 2\mathbb{N}$, by committing an error which is $o(1)$: more precisely, as $N \to \infty$, with $N \in 2\mathbb{N}$,

$$\mathcal{P}_{\delta,T_N}(\theta_1 > N-r, N-r \in \tau) = \mathcal{P}_{\delta,T_N}(\theta_1 > N-i-r, N-i-r \in \tau) + o(1).$$

Returning to (4.16) and replacing $\mathcal{P}_{\delta,T_N}(N \in \tau)$ by $\mathcal{P}_{\delta,T_N}(N - i \in \tau)$, we can write

$$\mathcal{P}_{\delta,T_N}(L'_N = j | N \in \tau) = \mathcal{P}_{\delta,T_N}(L'_{N-i} = j, \theta_j \leq N - \sqrt{N} | N - i \in \tau) + o(1)$$
$$= \mathcal{P}_{\delta,T_N}(L'_{N-i} = j | N - i \in \tau) + o(1),$$

where the second equality follows by (4.15). The proof is complete. $\square$

Let us return to (4.13). We write the event $\{\tau \cap (N - L, N] \neq \varnothing\}$ as a disjoint union

$$\{\tau \cap (N - L, N] \neq \varnothing\} = \bigcup_{i=0}^{L-1} \mathcal{A}_i,$$
(4.19)
$$\mathcal{A}_i := \{N - i \in \tau, N - k \notin \tau \text{ for } 0 \leq k < i\},$$

that is, $N - i$ is the last renewal epoch before $N$. We can then write the left-hand side of (4.13) as

(4.20)
$$\mathcal{P}_{\delta,T_N}(L'_N = j, \tau \cap (N - L, N] \neq \varnothing)$$
$$= \sum_{i=0}^{L-1} \mathcal{P}_{\delta,T_N}(L'_N = j | \mathcal{A}_i) \cdot \mathcal{P}_{\delta,T_N}(\mathcal{A}_i).$$

Note that $\mathcal{P}_{\delta,T_N}(L'_N = j | \mathcal{A}_i) = \mathcal{P}_{\delta,T_N}(L'_{N-i} = j | \mathcal{A}_i)$ because $L'_N = L'_{N-i}$ on the event $\mathcal{A}_i$. The next basic fact is that, by the renewal property, we have

$$\mathcal{P}_{\delta,T_N}(L'_{N-i} = j | \mathcal{A}_i) = \mathcal{P}_{\delta,T_N}(L'_{N-i} = j | N - i \in \tau)$$

because the event $\{L'_{N-i} = j\}$ depends only on $\tau \cap [0, N - i]$. Therefore, we can apply Lemma 11 and rewrite (4.20) as

$$\mathcal{P}_{\delta,T_N}(L'_N = j, \tau \cap (N - L, N] \neq \varnothing)$$
(4.21)
$$= \mathcal{P}_{\delta,T_N}(L'_N = j | N \in \tau) \left( \sum_{i=0}^{L-1} \mathcal{P}_{\delta,T_N}(\mathcal{A}_i) \right) + o(1)$$
$$= \mathcal{P}_{\delta,T_N}(L'_N = j | N \in \tau) \cdot \mathcal{P}_{\delta,T_N}(\tau \cap (N - L, N] \neq \varnothing) + o(1).$$

However, by (4.12), the term $\mathcal{P}_{\delta,T_N}(\tau \cap (N - L, N] \neq \varnothing)$ is as close to 1 as we wish, by taking $L$ large. Combining (4.13) with (4.21), this means that, for every $j \in \mathcal{N} \cup \{0\}$ and for $N$ sufficiently large, we have

$$\mathcal{P}_{\delta,T_N}(L'_N = j | N \in \tau) \in \left( e^{-t_{\delta,\varsigma}} \frac{(t_{\delta,\varsigma})^j}{j!} - 2\varepsilon, e^{-t_{\delta,\varsigma}} \frac{(t_{\delta,\varsigma})^j}{j!} + 2\varepsilon \right).$$

Since $\varepsilon$ is arbitrary, (4.11) is proved and the step is completed.

4.3. *Step 3.* In this last step, it remains to prove that for all $\varepsilon > 0$ and all $j \in \mathbb{N} \cup \{0\}$,

(4.22) $$\lim_{N \to \infty} \mathbf{P}_{N,\delta}^{T_N}\left(\frac{S_N}{T_N} \in [j-\varepsilon, j+\varepsilon]\right) = \mathbf{P}(S_\Gamma = j),$$

where $\Gamma$ is a random variable independent of the $\{S_i\}_{i \geq 0}$ and with a Poisson law of parameter $t_{\delta,\zeta}$.

Let $\varepsilon > 0$ and set

(4.23) $$V_\varepsilon(N) := \left|\mathbf{P}_{N,\delta}^{T_N}\left(\frac{S_N}{T_N} \in [j-\varepsilon, j+\varepsilon]\right) - \mathbf{P}(S_\Gamma = j)\right|.$$

Our goal is to prove that for all $\eta > 0$, we have $V_\varepsilon(N) \leq \eta$ when $N$ is large enough. We let $\mathcal{V}(N, l)$ be the set $\tau^{T_N} \cap [N-l, N]$ and it is useful to recall the result obtained in (3.27), that is, that there exists $\ell_0 = \ell_0(\eta)$ such that for every $N \geq \ell_0$, we have $\mathbf{P}_{N,\delta}^{T_N}(\mathcal{V}(N, \ell_0) = \varnothing) \leq \eta/4$. Therefore, with $N$ large enough, we obtain

$$V_\varepsilon(\delta) \leq \frac{\eta}{2} + \left|\mathbf{P}_{N,\delta}^{T_N}\left(\frac{S_N}{T_N} \in [j-\varepsilon, j+\varepsilon], \mathcal{V}(N, \ell_0) \neq \varnothing\right)\right.$$
$$\left. - \mathbf{P}(S_\Gamma = j)\mathbf{P}_{N,\delta}^{T_N}(\mathcal{V}(N, \ell_0) \neq \varnothing)\right|.$$

With some abuse of notation, we still denote by $\theta_j$ and $L'_N$ the variables on the $S$ space defined by (4.2) and (4.3) with $\tau_i$ replaced by $\tau_i^{T_N}$ and $\varepsilon_i$ by $\varepsilon_i^{T_N}$ (in particular, $L'_N := \#\{i \leq L_{N,T_N} : |\varepsilon_i^{T_N}| = 1\}$). Then, note that on the event $\mathcal{V}(N, \ell_0)$, we have $|S_N - S_{\theta_{L'_N}}| \leq \ell_0$. Moreover, for all $N \geq 1$, we have $S_{\theta_{L'_N}}/T_N \in \mathbb{Z}$. Therefore, assuming that $\varepsilon$ has been chosen small enough, we obtain, for $N$ large enough,

(4.24) $$\mathbf{P}_{N,\delta}^{T_N}\left(\frac{S_N}{T_N} \in [j-\varepsilon, j+\varepsilon], \mathcal{V}(N, \ell_0) \neq \varnothing\right)$$
$$= \mathbf{P}_{N,\delta}^{T_N}\left(\frac{S_{\theta_{L'_N}}}{T_N} = j, \mathcal{V}(N, \ell_0) \neq \varnothing\right).$$

We can rewrite the right-hand side of (4.24) by using, for $i \in \{0, \ldots, \ell_0\}$, the sets $\mathcal{A}_i$ introduced in (4.19). This gives

(4.25) $$\mathbf{P}_{N,\delta}^{T_N}\left(\frac{S_{\theta_{L'_N}}}{T_N} = j, \mathcal{V}(N, \ell_0) \neq \varnothing\right)$$
$$= \sum_{i=0}^{\ell_0} \mathbf{P}_{N,\delta}^{T_N}\left(\frac{S_{\theta_{L'_{N-i}}}}{T_N} = j \middle| \mathcal{A}_i\right) \mathbf{P}_{N,\delta}^{T_N}(\mathcal{A}_i).$$

At this stage, the Markov property and equation (2.13) give

$$\mathbf{P}_{N,\delta}^{T_N}(\cdot|\mathcal{A}_i) = \mathbf{P}_{N,\delta}^{T_N}(\cdot|N - i \in \tau) = \mathcal{P}_{\delta,T_N}(\cdot|N - i \in \tau).$$

Hence, we can rewrite (4.25) as

$$\begin{aligned}
(4.26) \quad & \mathbf{P}_{N,\delta}^{T_N}\left(\frac{S_{\theta_{L'_N}}}{T_N} = j, \mathcal{V}(N,\ell_0) \neq \varnothing\right) \\
& = \sum_{i=0}^{l_0} \mathcal{P}_{\delta,T_N}\left(\frac{S_{\theta_{L'_{N-i}}}}{T_N} = j \bigg| N - i \in \tau\right) \mathbf{P}_{N,\delta}^{T_N}(\mathcal{A}_i).
\end{aligned}$$

Thus, the proof of this step will be complete if we can show that, for all $i \in \{0, \ldots, \ell_0\}$,

$$\lim_{N \to \infty} \mathcal{P}_{\delta,T_N}\left(\frac{S_{\theta_{L'_{N-i}}}}{T_N} = j \bigg| N - i \in \tau\right) = \mathbf{P}(S_\Gamma = j).$$

This is proved once we show that, for all $(v,j) \in \mathbb{N} \cup \{0\} \times \mathbb{Z}$,

$$(4.27) \quad \lim_{N \to \infty} \mathcal{P}_{\delta,T_N}\left(L'_{N-i} = v, \frac{S_{\theta_v}}{T_N} = j \bigg| N - i \in \tau\right) = \mathbf{P}(\Gamma = v)\mathbf{P}(S_v = j).$$

We can rewrite the left-hand side of (4.27) as

$$(4.28) \quad \mathcal{P}_{\delta,T_N}\left(\frac{S_{\theta_v}}{T_N} = j \bigg| L'_{N-i} = v, N - i \in \tau\right) \cdot \mathcal{P}_{\delta,T_N}(L'_{N-i} = v | N - i \in \tau)$$

and it is easy to figure out that the process $(S_{\theta_n}/T_N)_{n \in \mathbb{N}}$ is just the symmetric simple random walk on $\mathbb{Z}$ and is independent of $(L'_{N-i}, \tau)$. Therefore, the first factor in (4.28) equals $\mathbf{P}(S_j = v)$ and then Lemma 11 and equation (4.11) are sufficient to complete the proof.

**5. Proof of Theorem 2(iii).** In this section, we prove part (iii) of Theorem 2. The parameter $\delta > 0$ is fixed throughout the section and the assumption is that the sequence $(T_N)_N$ is such that $T_N - \frac{\log N}{c_\delta} \to +\infty$ or, equivalently,

$$(5.1) \quad Q_{T_N}^1 \cdot N \longrightarrow 0 \quad (N \to \infty),$$

where we have used Lemma 7 and we recall the shorthand $Q_{T_N}^1 := Q_{T_N}^1(\phi(\delta, T_N))$ introduced in Section 3. The goal is to prove equation (1.15), that is, that the law of $S_N$ under $\mathbf{P}_{N,\delta}^{T_N}$ is tight.

In analogy with the previous sections, we start working under the law $\mathcal{P}_{\delta,T_N}$. We show that the polymer of length $N$ does not visit any interface other than the one located at $S = 0$, that is [recalling (4.2) and (4.3)], $L'_N = 0$. Note, in fact, that

$$\{L'_N \geq 1\} = \{\theta_1 \leq N\} = \bigcup_{i=1}^{L_N} \{|\varepsilon_i| = 1|\} \subseteq \bigcup_{i=1}^{N/2} \{|\varepsilon_i| = 1|\}$$

because, plainly, $L_N \leq N/2$; recall (3.2). Hence, the inclusion bound yields

(5.2)
$$\mathcal{P}_{\delta,T_N}(L'_N \geq 1) \leq \frac{N}{2} \cdot \mathcal{P}_{\delta,T_N}(|\varepsilon_1| = 1) = e^\delta N Q^1_{T_N}$$
$$\longrightarrow 0 \qquad (N \to \infty),$$

where we have used (4.5) and (5.1). With the same abuse of notation as in the previous section, we also denote by $L'_N$ the variable on the $S$ space defined by $L'_N := \#\{i \leq L_{N,T_N} : |\varepsilon_i^{T_N}| = 1\}$, so that applying (2.13), we get, as $N \to \infty$ with $N \in 2\mathbb{N}$,

(5.3)
$$\mathbf{P}^{T_N}_{N,\delta}(L'_N \geq 1 | S_N \in T_N\mathbb{Z}) = \mathcal{P}_{\delta,T_N}(L'_N \geq 1 | N \in \tau)$$
$$\leq \frac{\mathcal{P}_{\delta,T_N}(L'_N \geq 1)}{\mathcal{P}_{\delta,T_N}(N \in \tau)} \longrightarrow 0,$$

having applied (5.2) and Lemma 5.

Now, set $|S_n|^* := \max_{0 \leq k \leq n} |S_k|$ and observe that equation (5.3) can be rephrased as

$$\mathbf{P}^{T_N}_{N,\delta}(|S_N|^* \geq T_N | N \in \tau^{T_N}) \longrightarrow 0 \qquad (N \to \infty, N \in 2\mathbb{N}).$$

We want to remove the conditioning on $N \in \tau^{T_N}$. To this end, we let $\mu_N := \tau^{T_N}_{L_{N,T_N}}$ denote the location of the last point of $\tau^{T_N} \cap [0, N]$. Let us recall equation (3.27), which holds whenever $\delta > 0$ and which can hence be applied here: for every $\varepsilon > 0$, there exists $\ell_0$ such that $\mathbf{P}^{T_N}_{N,\delta}(\mu_N < N - \ell_0) < \varepsilon$ for every $N \in \mathbb{N}$. Therefore,

(5.4)
$$\left| \mathbf{P}^{T_N}_{N,\delta}(|S_N|^* \geq T_N) - \sum_{\ell=0}^{\ell_0} \mathbf{P}^{T_N}_{N,\delta}(|S_N|^* \geq T_N | \mu_N = N - \ell) \mathbf{P}^{T_N}_{N,\delta}(\mu_N = N - \ell) \right| \leq \varepsilon.$$

However, on the event $\{\mu_N = N - \ell\}$, we have $|S_N|^* \geq T_N$ if and only if $|S_{N-\ell}|^* \geq T_N$. Moreover, $\{|S_{N-\ell}|^* \geq T_N\} = \{L'_{N-\ell} \geq 1\}$, hence, using the Markov property and (2.13), for $\ell$ even, we get

$$\mathbf{P}^{T_N}_{N,\delta}(|S_N|^* \geq T_N | \mu_N = N - \ell)$$
$$= \mathbf{P}^{T_N}_{N,\delta}(L'_{N-\ell} \geq 1 | N - \ell \in \tau^{T_N})$$
$$= \mathcal{P}_{\delta,T_N}(L'_{N-\ell} \geq 1 | N - \ell \in \tau) \longrightarrow 0 \qquad (N \to \infty, N \in 2\mathbb{N}).$$

Then, equation (5.4) yields, for $N$ sufficiently large,

$$\mathbf{P}_{N,\delta}^{T_N}(|S_N|^* \geq T_N) \leq 2\varepsilon.$$

We can finally prove that $S_N$ is tight. Denoting by $\xi_N$ a quantity such that $|\xi_N| \leq 2\varepsilon$ for $N$ large, we have

$$\mathbf{P}_{N,\delta}^{T_N}(|S_N| \geq L) = \mathbf{P}_{N,\delta}^{T_N}(|S_N| \geq L, |S_N|^* < T_N) + \xi_N$$
$$\leq \mathbf{P}_{N,\delta}^{T_N}(\mu_N \leq N - L) + \xi_N,$$

where the inequality follows by the inclusion bound since $\{|S_N| \geq L, |S_N|^* < T_N\} \subseteq \{\mu_N \leq N - L\}$. Then, again by (3.27), if $L \geq \ell_0$, we have, for large $N$,

$$\mathbf{P}_{N,\delta}^{T_N}(|S_N| \geq L) \leq 3\varepsilon.$$

Since $\varepsilon > 0$ was arbitrary, it follows that

$$\lim_{L \to \infty} \sup_{N \in \mathbb{N}} \mathbf{P}_{N,\delta}^{T_N}(|S_N| \geq L) = 0,$$

hence (1.15) is proved and the proof is complete.

## APPENDIX A: COMPUTING $Q_T^I(\lambda)$

The computation of $Q_T^1(\lambda)$ and $Q_T^2(\lambda)$, defined in (2.4), is a classical problem; see [10], Chapter XIV. For completeness, here, we are going to derive explicit formulae for $Q_T^1(\lambda)$ and $Q_T^2(\lambda)$, using a simple martingale argument. We assume that $T \in \mathbb{N}$ (i.e., $T < \infty$).

For $\mu \in \mathbb{C}$ and $n \in \mathbb{N}$, we set

$$M_n := \frac{e^{\mu S_n}}{(\cosh \mu)^n}$$

and observe that the $\{M_n\}_{n \geq 0}$ under $\mathbf{P}$ is a $\mathbb{C}$-valued martingale (i.e., its real and imaginary parts are $\mathbb{R}$-valued martingales) with respect to the natural filtration of the simple random walk $\{S_i\}_i$. We will only be interested in the special cases when $\mu \in \mathbb{R}$ or $\mu \in (-\frac{\pi}{2}i, \frac{\pi}{2}i)$ so that, in any case, $\cosh \mu \in \mathbb{R}^+$ and therefore the expression $\log \cosh \mu$ is well defined with no need of further specifications.

We denote by $\mathbf{P}_1$ the law $\mathbf{P}(\cdot|S_1 = 1)$ and note that $\{M_n\}_{n \geq 1}$ is a martingale under $\mathbf{P}_1$. Moreover, both $\tau_1^T$ and $|\varepsilon_1^T|$ have the same law under $\mathbf{P}$ and $\mathbf{P}_1$. The optimal stopping theorem yields $\mathbf{E}_1(S_{\tau_1^T}) = \mathbf{E}_1(S_1)$, that is,

$$e^{\mu T}\mathbf{E}_1\left(\frac{1}{(\cosh \mu)^{\tau_1^T}}\mathbf{1}_{\{|\varepsilon_1^T|=1\}}\right) + \mathbf{E}_1\left(\frac{1}{(\cosh \mu)^{\tau_1^T}}\mathbf{1}_{\{|\varepsilon_1^T|=0\}}\right) = \frac{e^\mu}{\cosh \mu}.$$

Defining, for notational brevity, $Q_T^i := Q_T^i(\log \cosh \mu)$, we can rewrite this relation as

$$2e^{\mu T} Q_T^1 + Q_T^0 = \frac{e^\mu}{\cosh \mu}.$$

The analogous relation with $\mu$ replaced by $-\mu$ leads to the following couple of equations:

$$2\cosh(\mu T) Q_T^1 + Q_T^0 = 1,$$
$$2\sinh(\mu T) Q_T^1 = \tanh \mu,$$

which yield the solutions

(A.1)
$$Q_T^0(\log \cosh \mu) = 1 - \frac{\tanh(\mu)}{\tanh(\mu T)},$$
$$Q_T^1(\log \cosh \mu) = \frac{\tanh(\mu)}{2\sinh(\mu T)}$$

and for $Q_T(\cdot) := Q_T^0(\cdot) + 2Q_T^1(\cdot)$, we have

(A.2) $$Q_T(\log \cosh \mu) = 1 - \tanh(\mu) \cdot \frac{\cosh(\mu T) - 1}{\sinh(\mu T)}.$$

Setting $\lambda = \log \cosh \mu$, that is, $\mu = \lambda + \log(1 + \sqrt{1 - e^{-2\lambda}})$, we finally obtain

(A.3)
$$Q_T^0(\lambda) = 1 - \sqrt{1 - e^{-2\lambda}} \cdot \frac{(1 + \sqrt{1 - e^{-2\lambda}})^T + (1 - \sqrt{1 - e^{-2\lambda}})^T}{(1 + \sqrt{1 - e^{-2\lambda}})^T - (1 - \sqrt{1 - e^{-2\lambda}})^T},$$
$$Q_T^1(\lambda) = \frac{\sqrt{1 - e^{-2\lambda}} \cdot e^{-\lambda T}}{(1 + \sqrt{1 - e^{-2\lambda}})^T - (1 - \sqrt{1 - e^{-2\lambda}})^T}$$

and, therefore,

(A.4)
$$Q_T(\lambda) = 1 - \sqrt{1 - e^{-2\lambda}}$$
$$\cdot \frac{(1 + \sqrt{1 - e^{-2\lambda}})^T + (1 - \sqrt{1 - e^{-2\lambda}})^T - 2e^{-\lambda T}}{(1 + \sqrt{1 - e^{-2\lambda}})^T - (1 - \sqrt{1 - e^{-2\lambda}})^T}.$$

Note that when $\lambda < 0$, we have $\mu = \lambda + \log(1 + \sqrt{1 - e^{-2\lambda}}) = i \arctan \sqrt{e^{-2\lambda} - 1}$, hence we can write, more explicitly,

(A.5)
$$Q_T^0(\lambda) = 1 - \frac{\sqrt{e^{-2\lambda} - 1}}{\tan(T \arctan \sqrt{e^{-2\lambda} - 1})},$$
$$Q_T^1(\lambda) = \frac{\sqrt{e^{-2\lambda} - 1}}{2\sin(T \arctan \sqrt{e^{-2\lambda} - 1})},$$
$$Q_T(\lambda) = 1 + \sqrt{e^{-2\lambda} - 1} \cdot \frac{1 - \cos(T \arctan \sqrt{e^{-2\lambda} - 1})}{\sin(T \arctan \sqrt{e^{-2\lambda} - 1})}.$$

A POLYMER IN A MULTI-INTERFACE MEDIUM 37

Of course, these formulae break down if $|\lambda|$ is too large. This happens at the first negative zero of the denominator $\lambda = \lambda_0^T$, where $(T \arctan \sqrt{e^{-2\lambda_0^T} - 1}) = \pi$, that is,

$$\text{(A.6)} \qquad \lambda_0^T := -\frac{1}{2} \log\left(1 + \left(\tan \frac{\pi}{T}\right)^2\right).$$

Note that as $\lambda \downarrow \lambda_o^T$, both $Q_T^0(\lambda)$ and $Q_T^1(\lambda)$ diverge (they have a pole). Also, note that taking the limit $\lambda \to 0$ in (A.3) or (A.5), we get

$$Q_T^0(0) = 1 - \frac{1}{T}, \qquad Q_T^1(0) = \frac{1}{2T}.$$

We conclude by noting that the probabilities $q^j(n)$ introduced in (2.4) can also be given explicit formulae. More precisely, by equation (5.8) in Chapter XIV of [10], we have, for all $n \geq 2$,

$$\text{(A.7)} \qquad \begin{aligned} q_T^0(n) &= \left(\frac{2}{T} \sum_{\nu=1}^{\lfloor (T-1)/2 \rfloor} \cos^{n-2}\left(\frac{\pi \nu}{T}\right) \sin^2\left(\frac{\pi \nu}{T}\right)\right) \cdot \mathbf{1}_{\{n \text{ is even}\}}, \\ q_T^1(n) &= \left(\frac{1}{T} \sum_{\nu=1}^{\lfloor (T-1)/2 \rfloor} (-1)^{\nu+1} \cos^{n-2}\left(\frac{\pi \nu}{T}\right) \sin^2\left(\frac{\pi \nu}{T}\right)\right) \cdot \mathbf{1}_{\{n-T \text{ is even}\}}. \end{aligned}$$

DIPARTIMENTO DI MATEMATICA PURA E APPLICATA
UNIVERSITÀ DEGLI STUDI DI PADOVA
VIA TRIESTE 63
35121 PADOVA
ITALY
E-MAIL: francesco.caravenna@math.unipd.it

EURANDOM
P.O. BOX 513
5600 MB EINDHOVEN
THE NETHERLANDS
E-MAIL: petrelis@eurandom.tue.nl